\documentclass[oneside,english]{amsart}

\usepackage[latin1]{inputenc} 
\usepackage{setspace} 
\usepackage{amssymb}
\usepackage[]{epsf,epsfig,amsmath,amssymb,amsfonts,latexsym}

\usepackage{amsthm}
\usepackage{mathtools}
\usepackage{lmodern}
\usepackage{mathabx}
\usepackage{stmaryrd}
\usepackage{mathrsfs}
\usepackage{wasysym}
\usepackage{bbm}		%
\usepackage[inline]{enumitem}
\usepackage{nccmath}	%
\usepackage{braket}		%
\usepackage[normalem]{ulem}		%
\usepackage{nicefrac}
\usepackage{cancel}
\usepackage{xcolor}
\usepackage[a4paper, left=2.7cm, right=2.7cm, top=3.5cm, bottom=2cm]{geometry}
\usepackage[pdftex=true,hyperindex=true,pdfborder={0 0 0}]{hyperref}
\usepackage{doi}
\usepackage[numbers]{natbib}
\usepackage{cleveref}
\usepackage{leftidx} %

\theoremstyle{plain}
\newtheorem{theorem}{Theorem}[section]
\newtheorem{lemma}[theorem]{Lemma}
\newtheorem{corollary}[theorem]{Corollary}
\newtheorem{proposition}[theorem]{Proposition}

\newtheorem{definition}[theorem]{Definition}

\newtheorem{maintheorem}{Theorem}
\newtheorem*{pirillotheorem*}{Pirillo's Theorem}

\theoremstyle{definition}
\newtheorem{remark}[theorem]{Remark}

\newcommand{\ZZ}{\mathbb{Z}}			%
\newcommand{\NN}{\mathbb{N}}			%
\newcommand{\RR}{\mathbb{R}}			%
\newcommand{\QQ}{\mathbb{Q}}			%

\newcommand{\symb}[1]{\mathtt{#1}}		%
\newcommand{\isdef}{:=}			%
\DeclarePairedDelimiter\norm{\lVert}{\rVert}	%

\newcommand{\indicator}[1]{\mathbbm{1}_{#1}}					%

\newcommand{\Lcal}{\mathcal{L}}

\newcommand{\occ}{\mathsf{occ}}
\newcommand{\zeropoint}{\boldsymbol{.}}

\newcommand{\define}[1]{\textbf{#1}}

\title{%
    A characterization of Sturmian sequences by
    indistinguishable asymptotic pairs
}

\author{Sebasti\'an Barbieri, S\'ebastien Labb\'e and \v{S}t\v{e}p\'{a}n Starosta}

\newcommand{\Addresses}{{
		\bigskip

		\hskip-\parindent   S.~Barbieri, \textsc{DMCC, Universidad de Santiago de Chile,
			Las Sophoras 173. Estaci\'on Central. Santiago. Chile.}\par\nopagebreak
		\textit{E-mail address}: \texttt{sebastian.barbieri@usach.cl}
		
		\medskip
		
		\hskip-\parindent   S.~Labb\'{e}, \textsc{LaBRI, Universit\'{e} de Bordeaux,
			351, cours de la Lib\'{e}ration, F-33405, Talence, France.}\par\nopagebreak
		\textit{E-mail address}: \texttt{sebastien.labbe@labri.fr}
		
		\medskip
		
		\hskip-\parindent   \v{S}.~Starosta, \textsc{Faculty of Information Technology, Czech Technical University in Prague, Th\'akurova 9,
160 00 Praha 6, Czech Republic}\par\nopagebreak
		\textit{E-mail address}: \texttt{stepan.starosta@fit.cvut.cz}
	}}
	
	\date{}

\begin{document}
	
		\begin{abstract}
            We give a new characterization of biinfinite Sturmian sequences in terms of indistinguishable asymptotic pairs. Two asymptotic sequences on a full $\mathbb{Z}$-shift are indistinguishable if the sets of occurrences of every pattern in each sequence coincide up to a finitely supported permutation. This characterization can be seen as an extension to biinfinite sequences of Pirillo's theorem which characterizes Christoffel words. Furthermore, we provide a full characterization of indistinguishable asymptotic pairs on arbitrary alphabets using substitutions and biinfinite characteristic Sturmian sequences. The proof is based on the well-known notion of derived sequences.
		\end{abstract}

		\maketitle

		\medskip	
		
        \smallskip
		
		\textit{Secondary:}
		37C29, %

\section{Introduction}

\begin{figure}[h]
    \begin{center}
        \includegraphics{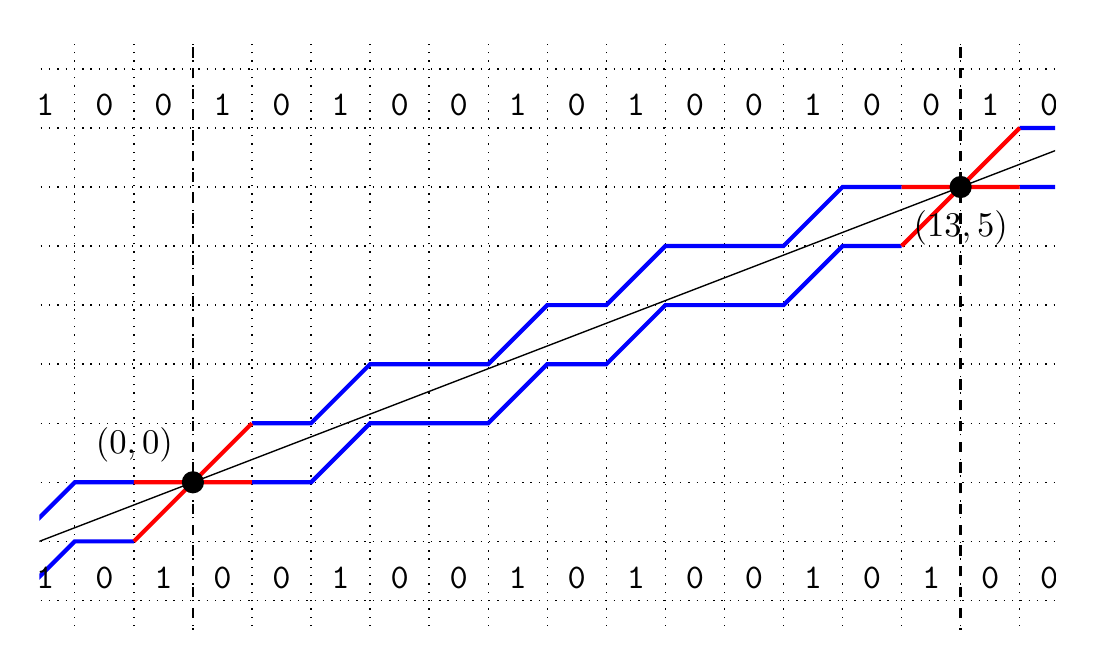}
    \end{center}
    \caption{The lower and upper sequences $c_\alpha$ and $c'_\alpha$ when
    $\alpha=5/13$ are periodic.}
    \label{fig:sturmian}
    \end{figure}

    Let $\alpha\in [0,1]$ and consider
    the lower and upper sequences $c_{\alpha}$ and $c'_{\alpha}$ given respectively by
	\[
	\begin{array}{rccl}
	c_{\alpha}:&\ZZ & \to & \{\symb{0},\symb{1}\}\\
	&n & \mapsto &\lfloor\alpha(n+1)\rfloor-\lfloor\alpha n\rfloor
	\end{array}
	\quad
	\text{ and }
	\quad
	\begin{array}{rccl}
	c'_{\alpha}:&\ZZ & \to & \{\symb{0},\symb{1}\}\\
	&n & \mapsto &\lceil\alpha(n+1)\rceil-\lceil\alpha n\rceil.
	\end{array}
	\]
    When $\alpha$ is rational, the sequences
    $c_{\alpha}$ and $c'_{\alpha}$ are periodic 
    and their period corresponds to Christoffel words \cite{MR2464862},
    see Figure~\ref{fig:sturmian}.
    More precisely, the shortest periodic pattern and smallest for the
    lexicographic order of $c_\alpha$ is the lower Christoffel word of slope
    $p/q$
    where $p$ and $q$ are nonnegative coprime integers such that
    $\alpha=p/(p+q)$.
    For example, when $\alpha=5/13$,
    the lower sequence $c_\alpha$ has period $\symb{0010010100101}$
    which is the lower Christoffel word of slope $5/8$
    and
    the upper sequence $c'_\alpha$ has period $\symb{1010010100100}$
    which is the upper Christoffel word of slope $5/8$.
    When $\alpha$ is irrational, then $c_\alpha$ and $c'_\alpha$ are 
    not periodic.
    The restrictions of $c_\alpha$ and $c'_\alpha$ to $\ZZ_{\geq1}$ are equal
    and correspond to the
    well-known one-sided characteristic Sturmian sequence
    of slope $\alpha$~\cite{MR0000745}.
    In this work, we consider biinfinite sequences as opposed to
    one-sided sequences.
    Over the domain $\ZZ$, we say that $c_\alpha$ and $c'_\alpha$ are
    respectively the \textbf{lower} and \textbf{upper characteristic Sturmian
    sequences of slope $\alpha$} whenever $\alpha$ is irrational.

Sturmian sequences have many equivalent definitions, for example, in terms of
aperiodic balanced sequences \cite{MR0000745}, 
irrational rotations \cite{MR1970391,MR1905123},
factor complexity \cite{MR0322838}
or return words \cite{MR1808196}.
On the other hand, Christoffel words also have many equivalent definitions,
including 14 characterizations listed in \cite{MR2681739}, see also
\cite{MR1453849,MR2464862}.
A recent book \cite{MR3887697} gathers exhaustively the combinatorial properties of
Christoffel words and uses them to prove two important theorems of Markoff for
Diophantine approximations and quadratic forms \cite{MR1510073}.

    In this work, we study a surprising connection between Sturmian sequences
    and asymptotic pairs satisfying a natural combinatorial property which
    originates in thermodynamical formalism. This property characterizes asymptotic pairs which induce the trivial linear functional on a space of continuous and shift-invariant cocycles on the asymptotic relation of the full $\ZZ$-shift. See Section 3 of~\cite{Barbieri_Gomez_Marcus_Meyerovitch_Taati_2020} for further details.

    Concretely, given a finite set $\Sigma$, we consider the space of
    sequences $\Sigma^{\ZZ} = \{ x \colon \ZZ \to \Sigma\}$ endowed with
    the prodiscrete topology and the shift action $\ZZ
    \overset{\sigma}{\curvearrowright} \Sigma^{\ZZ}$. In this setting, two
    sequences $x,y \in \Sigma^{\ZZ}$ are \define{asymptotic} if $x$ and
    $y$ differ in finitely many positions of $\ZZ$. The finite set $F = \{ n \in
    \ZZ : x_n \neq y_n\}$ is called the \define{difference set} of
    $(x,y)$.

    Given two asymptotic sequences $x,y \in \Sigma^{\ZZ}$ with the difference set $F$, we want
    to compare the number of occurrences of a fixed pattern. 
    As $x$ and $y$ are asymptotic, occurrences of patterns whose support do not
    intersect $F$ 
    are the same, so we only need to consider
    the occurrences of patterns that appear intersecting $F$.
    As an example, we can take a fixed symbol $a \in \Sigma$ and
    define $\Delta_a(x,y)$ as the number of positions $n \in F$ such that
    $y_n=a$ minus the number of positions $n\in F$ such that $x_n=a$. 
    As $F$ is finite, this value is well defined. 
    More generally, for any given pattern $p \colon S
    \to \Sigma$ where $S$ is a finite subset of $\ZZ$, we can consider the
    difference $\Delta_p(x,y)$ of the number of occurrences of $p$ in $y$
    intersecting $F$ minus the number of occurrences of $p$ in $x$ intersecting
    $F$.

    We say that $(x,y)$ is
    an \define{indistinguishable asymptotic pair} 
    if $(x,y)$ is asymptotic and
    $\Delta_p(x,y) =0$ for every pattern $p$. 
    A trivial example of an indistinguishable
    asymptotic pair is $(x,x)$ for any $x \in \Sigma^{\ZZ}$. Another simple
    example is $x,y \in \{\symb{0},\symb{1}\}^{\ZZ}$ where $x$ is equal to
    $\symb{1}$ at the origin, and $\symb{0}$ everywhere else, and $y$ is equal
    to $\symb{1}$ at some nonzero $n \in \ZZ$ and $\symb{0}$ everywhere else.
    Note that in both of these examples $x$ and $y$ lie on the same orbit of
    $\ZZ \overset{\sigma}{\curvearrowright} \Sigma^{\ZZ}$. 
    
    In~\cite{Barbieri_Gomez_Marcus_Meyerovitch_Taati_2020} the authors define the following norm on asymptotic sequences of $\Sigma^{\ZZ}$ \[ \norm{(x,y)}^{*}_{\mathsf{NS}} = \sup_{\substack{S \subseteq \ZZ\\S \mbox{ finite}}} \frac{1}{|S|}\sum_{p \in \Sigma^S}|\Delta_p(x,y)|.  \]

	Every asymptotic pair induces an evaluation map on the space of continuous cocycles on the equivalence relation of asymptotic pairs. The authors show that this norm coincides with the dual norm in the space of linear functionals on the space of continuous cocycles. In other words, the asymptotic pairs which induce the trivial linear functional are precisely the indistinguishable pairs. In this article, we provide a full characterization of which asymptotic pairs induce the trivial linear functional.
	
    Using the notion of indistinguishability, 
    we provide a characterization of
    the lower and the upper characteristic Sturmian
    sequences.

	\begin{maintheorem}\label{thm:sturmian_characterization}
		Let $x,y\in\{\symb{0},\symb{1}\}^\ZZ$ and assume that $x$ is recurrent.
		The pair $(x,y)$ is an indistinguishable asymptotic pair
		with difference set $F=\{-1,0\}$ such that $x_{-1}x_{0}=\symb{10}$
		and $y_{-1}y_{0}=\symb{01}$ if and only if
            there exists $\alpha\in[0,1]\setminus\QQ$ such that
			$x={c}_{\alpha}$ and $y={c}'_{\alpha}$ are the lower and upper characteristic
			Sturmian sequences of slope $\alpha$.
	\end{maintheorem}

    \Cref{thm:sturmian_characterization} is proved in \Cref{sec:sturmian}.
    The auxiliary claims in this section bring a new understanding of Sturmian sequences. 
    In particular,
    for all $n\in\NN$ the words
    $x_{-n}x_{-n+1}\cdots x_{n-1}$
    and
    $y_{-n}y_{-n+1}\cdots y_{n-1}$
    of length $2n$
    are optimal representations of the language of $c_\alpha$ as
    both contain exactly one occurrence of every factor of length $n$,
    see \Cref{cor:two-compact-representation-of-the-language}.
    Removing the hypothesis that $x$ is recurrent, we
    obtain a unifying description of the lower and upper characteristic Sturmian sequences and
    their limits as their slope tends towards a rational value.

	\begin{maintheorem}\label{thm:sturmian_characterization_Z1}
		Let $x,y\in\{\symb{0},\symb{1}\}^\ZZ$.
		The pair $(x,y)$ is an indistinguishable asymptotic pair
		with difference set $F=\{-1,0\}$ such that $x_{-1}x_{0}=\symb{10}$
		and $y_{-1}y_{0}=\symb{01}$ if and only if
        there exists $(\alpha_n)_{n \in \NN}$ with $\alpha_n \in [0,1]\setminus \QQ$ such that \[
		x=\lim_{n \to \infty}c_{\alpha_n} \quad \mbox{ and }\quad y=\lim_{n \to \infty}c'_{\alpha_n}.\]
	\end{maintheorem}
	
    In the case where $x$ is not recurrent,
    then $x$ and $y$ lie on the same orbit and there exist coprime integers
	$p,q\in\ZZ_{\geq0}$ such that
	$(x,y)$ is the limit of asymptotic pairs formed by the lower and upper characteristic Sturmian
    sequences of slope $\alpha_n$ as $\alpha_n$ converges toward the rational
    slope $p/(p+q)\in[0,1]\cap\QQ$ either
    from above or from below,
    see \Cref{thm:sturmian_characterization_non_recurrent}.
    Limits of the lower and upper characteristic Sturmian sequences as their slope tends to a rational number are expressed in terms of Christoffel words, see
    \Cref{lem:limit-of-sturmian-as-christoffel}.
	The proof of \Cref{thm:sturmian_characterization_Z1} follows from
    \Cref{thm:sturmian_characterization} and
    \Cref{thm:sturmian_characterization_non_recurrent}
    and is proved in~\Cref{sec:limit-sturmian}. 

    \Cref{thm:sturmian_characterization} and \Cref{thm:sturmian_characterization_Z1}
    are related to a famous theorem of Pirillo \cite{MR1854493}
    which provides a characterization of Christoffel words
    of slope $p/q$ where
    $p$ and $q$ are positive coprime integers.
If $p$ and $q$ are nonzero, the lower Christoffel word of slope $p/q$ starts with
letter $\symb{0}$ and ends with letter $\symb{1}$, so it can be written as $\symb{0}m\symb{1}$
for some finite word $m\in\{\symb{0},\symb{1}\}^*$
and the corresponding upper Christoffel word is $\symb{1}m\symb{0}$.
Pirillo gave the following elegant characterization of Christoffel words.
Recall that two words $w,w'\in\{\symb{0},\symb{1}\}^*$ are \define{conjugate} if there exists
$u,v\in\{\symb{0},\symb{1}\}^*$ such that $w=uv$ and $w'=vu$.

\begin{pirillotheorem*}[{\cite{MR1854493}}]
\label{thm:pirillo}
    The word $\symb{0}m\symb{1}\in\{\symb{0},\symb{1}\}^*$ is a lower Christoffel word if
    and only if $\symb{0}m\symb{1}$ and $\symb{1}m\symb{0}$ are conjugate.
\end{pirillotheorem*}

    \begin{figure}[h]
    \begin{center}
        \includegraphics{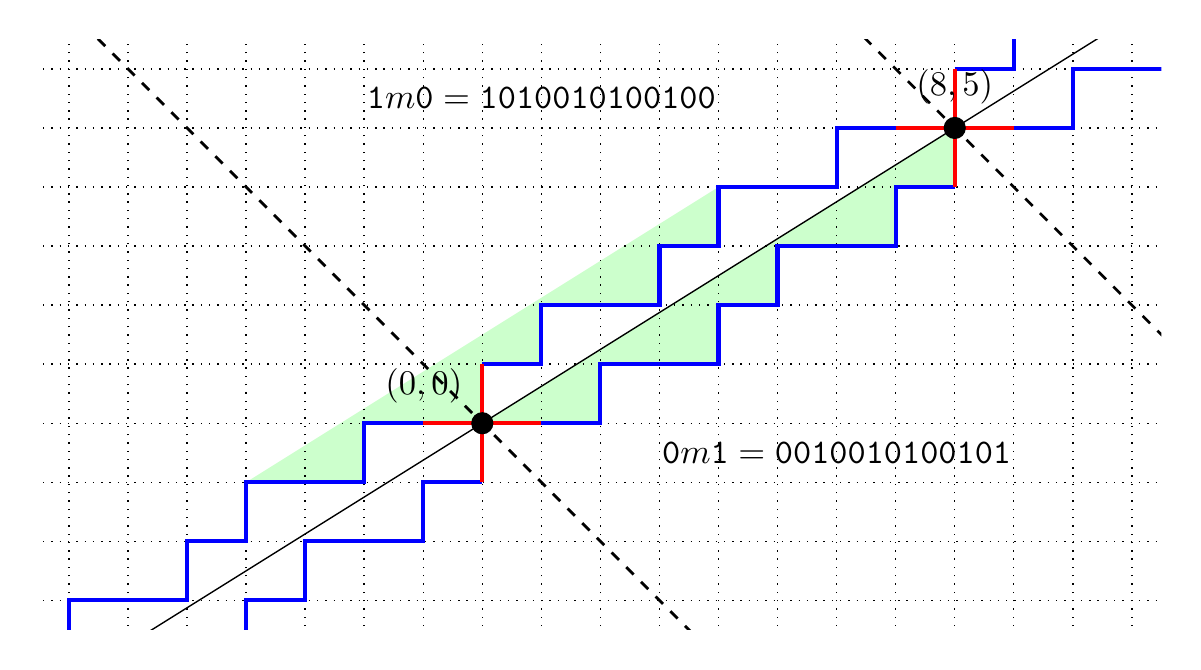}
    \end{center}
    \caption{Pirillo's theorem characterizes
    Christoffel words: the lower Christoffel word 
    $\symb{0}m\symb{1}\in\{\symb{0},\symb{1}\}^*$
    is conjugate to the upper Christoffel word $\symb{1}m\symb{0}$.}
    \label{fig:Pirillo}
    \end{figure}

Pirillo's theorem is illustrated in Figure~\ref{fig:Pirillo}.
We observe that the conjugacy of 
$\symb{0}m\symb{1}$
into $\symb{1}m\symb{0}$ is done via their factorization into a product of two
palindromes:
    $\symb{0}m\symb{1}=\symb{00100}\cdot\symb{10100101}$ and
    $\symb{1}m\symb{0}=\symb{10100101}\cdot\symb{00100}$.
The factorization of $\symb{0}m\symb{1}$ as a product of two palindromes and
the fact that the central word $m$ is a palindrome \cite[Prop.~4.2]{MR2464862} is also a
characterization of Christoffel words, see \cite{MR1311214} and \cite[Theorem
12.2.10]{MR3887697}.

    Pirillo's theorem can be restated for biinfinite sequences as follows:
    \emph{$c_\alpha$ is the lower sequence associated to the
        rational slope $\alpha=p/(p+q)$
        for some coprime nonnegative integers $p,q$
        if and only if $c_\alpha$ is a shift of $c'_\alpha$}.
    It is natural to ask if there is an analogous statement which holds as we take the limit
    $\frac{p}{p+q}\to\alpha$ for some irrational $\alpha\in[0,1]\setminus\QQ$. In this light,~\Cref{thm:sturmian_characterization_Z1} 
    can be considered as the extension of Pirillo's theorem to
    aperiodic biinfinite sequences where
    the notion of conjugacy of words is replaced
    by the notion of indistinguishability of an asymptotic pair.
    This seems to be the correct approach
    since other alternatives (e.g., having the same language, see~\Cref{rem:toeplitz}) fail.

    The next result provides a full characterization of non-trivial indistinguishable asymptotic pairs for $\ZZ$ which does not depend upon the form of the difference set or the alphabet. More precisely, we show that every indistinguishable asymptotic pair can be obtained from limits of pairs of lower and upper characteristic Sturmian sequences by means of shifts and substitutions.
    
    Given finite sets $\Sigma,\Gamma$, a map $\varphi\colon \Sigma \to \Gamma^+$ which replaces symbols of $\Sigma$ by nonempty words on $\Gamma$ is called a substitution. This map is naturally extended by concatenation to a continuous map $\varphi \colon \Sigma^{\ZZ} \to \Gamma^{\ZZ}$.
    
    \begin{maintheorem}\label{thm:characterization_Z}
    	Let $\Sigma$ be a finite alphabet and $x,y \in \Sigma^{\ZZ}$ a non-trivial asymptotic pair.
        Then $x,y$ is indistinguishable if and only if either
    	\begin{itemize}
    		\item $x$ is recurrent and
    		there exists $\alpha\in[0,1]\setminus\QQ$,
    		a substitution 
    		$\varphi\colon \{\symb{0},\symb{1}\} \to \Sigma^+$ 
    		and an integer $m\in\ZZ$
    		such that 
    		\[  
    		\{x,y  \} = \{\sigma^m\varphi(\sigma({c}_{\alpha})), 
                         \sigma^m\varphi(\sigma({c}'_{\alpha}))  \},
    		\]
    		\item $x$ is not recurrent and
    		there exists a substitution 
    		$\varphi\colon \{\symb{0},\symb{1}\} \to \Sigma^+$  and an integer $m\in\ZZ$
    		such that \[\{x,y  \} = \{\sigma^m\varphi(\leftidx^{\infty}\symb{0}\zeropoint\symb{10}^{\infty}), \sigma^m\varphi(\leftidx^{\infty}\symb{0}\zeropoint\symb{010}^{\infty})  \}.\]
    	\end{itemize}
    \end{maintheorem}

    This means that every indistinguishable asymptotic pair in $\ZZ$ consists either of (1) two sequences in the same orbit, which are shifts of a sequence of the form $\leftidx^{\infty}v.uv^{\infty}$ for some $u,v \in \Sigma^+$, or (2) two sequences which, up to translation, can be obtained through a substitution from a pair of lower and upper characteristic Sturmian sequences. In simpler terms, all non-trivial examples of one-dimensional indistinguishable asymptotic pairs arise from irrational circle rotations. The proof of~\Cref{thm:characterization_Z} is given in~\Cref{sec:resultsZ}.
    It is based on the well-known notions of return words and derived sequences
    \cite{MR1489074}.

	\textbf{Acknowledgments}: 
    The first two authors were supported by the Agence Nationale de la Recherche through the
    projects CODYS (ANR-18-CE40-0007) and CoCoGro (ANR-16-CE40-0005). S. Barbieri was also supported by the FONDECYT grant 11200037. \v{S}. Starosta acknowledges the support of the OP VVV MEYS funded project
	CZ.02.1.01/0.0/0.0/16\_019/0000765  ``Research Center for Informatics''.
    This work originated from a visit of the first two authors to Prague in
    October 2019 supported by PHC Barrande, a France-Czech Republic bilateral
    funding and grant no. 7AMB18FR048 of MEYS of Czech Republic.
	
	\section{Preliminaries}\label{sec:preliminaries}

	Let $\NN$ denote the set of nonnegative integers. Intervals consisting of integers will be written using the notation $\llbracket n,m \rrbracket = [n,m]\cap \ZZ$, for $n,m \in \ZZ$.
	
	Let $\Sigma$ be a finite set to which we refer as an \define{alphabet}. An element $x \in \Sigma^{\ZZ} = \{x \colon \ZZ \to \Sigma  \}$ is called a biinfinite \define{sequence}. We shall often omit the word ``biinfinite'' and use the word ``one-sided'' to refer to functions with domain $\ZZ_{\geq1}$. For $n \in \ZZ$, we write $x_n$ to denote the value $x(n)$. The set $\Sigma^{\ZZ}$ of all sequences is endowed with the prodiscrete topology, which is generated by the metric \[ d(x,y) = 2^{- \inf\{|n|\ :\ n \in \ZZ  \mbox{ and } x_n \neq y_n \} }.   \]

    The \define{shift} is the map $\sigma \colon \Sigma^{\ZZ}\to \Sigma^{\ZZ}$ where
    \[ \left(   \sigma(x) \right)_m  =  x_{m+1} \quad \mbox{ for every } m \in \ZZ \mbox{ and } x \in \Sigma^{\ZZ}.  \]

	Let us represent pictorially a sequence $x \in \Sigma^\ZZ$ by marking the position of the zero coordinate with a point as follows:
	\[ x = \dots x_{-5}x_{-4}x_{-3}x_{-2}x_{-1} \zeropoint x_{0}x_{1}x_{2}x_{3}x_{4}\dots  \]
	Given words $u,w \in \Sigma^+$ and $y,z \in \Sigma^*$ we shall use the notation \[ x = \leftidx^{\infty}uy\zeropoint zw^{\infty} \in \Sigma^{\ZZ}  \]
	to indicate that the sequence consists of repeated concatenations of $u$ to the left of $y$, and of repeated concatenations of $w$ to the right of $z$.
	
	\begin{definition}
		We say that two sequences $x,y$ are \define{asymptotic}, or that $(x,y)$ is an asymptotic pair, if the set $F = \{ n \in \ZZ \colon x_n \neq y_n\}$ is finite. $F$ is called the \define{difference set} of $(x,y)$. If $x=y$ we say that the asymptotic pair is \define{trivial}.
	\end{definition}
	
	Equivalently, $x,y$ are asymptotic if for every sequence $\{n_i\}_{i \in \NN}$ of elements of $\ZZ$ such that $|n_i| \to \infty$, the distance $d(\sigma^{n_i}(x),\sigma^{n_i}(y))$ converges to zero.
	
	For finite $S \subseteq \ZZ$, a function $p\colon S \to \Sigma$ is called a \define{pattern} and the set $S$ is its \define{support}. Given a pattern $p \in \Sigma^S$, the \define{cylinder} centered at $p$ is $[p] = \{ x \in \Sigma^{\ZZ} \colon x|_S = p \}$. We say a pattern $p$ \define{appears} in $x \in \Sigma^{\ZZ}$ if there exists $n \in \ZZ$ such that $\sigma^n(x) \in [p]$. Let us also denote by $\occ_p(x) = \{n \in \ZZ \colon \sigma^n(x) \in [p]\}$ the set of \define{occurrences} of $p$ in $x \in \Sigma^{\ZZ}$.

    Given a pattern $p$ and sequences $x,y \in \Sigma^{\ZZ}$, we want to define a number which counts the difference between the occurrences of $p$ in $y$ compared to the occurrences of $p$ in $x$. Naively, if $\indicator{[p]}$ is the indicator function of $[p]$, we would like to sum over all integers $n$
    the difference $\indicator{[p]}(\sigma^n(y))-\indicator{[p]}(\sigma^n(x))$.
    For arbitrary $x,y \in \Sigma^{\ZZ}$ this sum is not well defined. However, it can be given meaning if $x,y$ are asymptotic. Indeed, if $F$ denotes the difference set of $x,y$ and $S$ denotes the support of $p$, then $F-S = \left \{ f-s : f \in F, s \in S \right \}$ is the set of all integers for which there exists an $s \in S$ such that $n+s \in F$. In consequence, we have that if $n \in \ZZ \setminus (F-S)$ then for every $s \in S$ we have that $n+s \notin F$ and thus $\sigma^n(x)_s = \sigma^n(y)_s$, which implies in turn that $\indicator{[p]}(\sigma^n(x)) = \indicator{[p]}(\sigma^n(y))$. This motivates the following definition.
	
	\begin{definition}
		Let $p$ be a pattern with finite support $S \subseteq \ZZ$ and $x,y \in \Sigma^{\ZZ}$ be asymptotic sequences with difference set $F$. The \define{discrepancy} of the pattern $p$ associated to the pair $(x,y)$ is given by \[ \Delta_p (x,y) = \sum_{n \in F-S } \indicator{[p]}(\sigma^n(y))-\indicator{[p]}(\sigma^n(x)).\]
	\end{definition}

    For example, 
    the discrepancy of the pattern $p=abcabc$ in the sequences
    \begin{align*}
        x &= \cdots bcabcbc\underline{abcabc}bc\underline{abc\zeropoint} \framebox{$\underline{abc}$}bcabcbc\underline{abcabc}bcabcbc\cdots,\\
        y &= \cdots bcabcbc\underline{abcabc}bcabc\zeropoint \framebox{$bc\underline{a}$}\underline{bcabc}bc\underline{abcabc}bcabcbc\cdots
    \end{align*}
    is $\Delta_{abcabc}(x,y) = 1-1=0$, because both $x$ and $y$ contain exactly
    one occurrence of the pattern $p=abcabc$ intersecting the difference set
    $F=\{0,1,2\}$.

	\begin{definition}
        We say that an asymptotic pair $(x,y)$ is an \define{indistinguishable
        asymptotic pair} if the discrepancy of every pattern $p$ of finite
        support is $\Delta_p(x,y) = 0$.
	\end{definition}

    A related notion is the one of \emph{local indistinguishability} given in
    \cite[\S~5.1]{baake_aperiodic_2013} which corresponds symbolically to having two
    sequences with the same language (as defined in the next section). 
    This notion applies to a more general context as it can be defined for pairs which are not asymptotic.
    An indistinguishable asymptotic
    pair is locally indistinguishable but the converse is not true even for
    asymptotic pairs as explained in Remark~\ref{rem:toeplitz}.

	Whenever $x,y$ are asymptotic, for every pattern $p$ the sets $\occ_p(x)$ and $\occ_p(y)$ are asymptotic when seen as sequences in $\{\symb{0},\symb{1}\}^{\ZZ}$. In consequence, the condition $\Delta_p(x,y)=0$ is equivalent to having $\occ_p(x)$ and $\occ_p(y)$ coincide up to a finitely supported permutation of $\ZZ$. More precisely, we have:\[ \#(\occ_p(x) \cap (F-S)) = \#(\occ_p(y) \cap (F-S)).   \]
	
	We are interested in understanding which asymptotic pairs are indistinguishable. In order to avoid simple cases, we will restrict our search to asymptotic pairs which are non-trivial. Notice that non-trivial indistinguishable asymptotic pairs may consist of sequences that lie in the same orbit. For example,
		\begin{align*}
		x = 
		\leftidx^{\infty}\symb{0}\symb{0}\symb{0}\symb{0} \framebox{$\symb{0}\symb{0}\zeropoint\symb{1}\symb{1}$} \symb{0}\symb{0}\symb{0}\symb{0}^{\infty} 
            \qquad \text{ and }\qquad
		y = 
		\leftidx^{\infty}\symb{0}\symb{0}\symb{0}\symb{0} \framebox{$\symb{1}\symb{1}\zeropoint\symb{0}\symb{0}$} \symb{0}\symb{0}\symb{0}\symb{0}^{\infty} 
\;
		\end{align*}
        is an indistinguishable asymptotic pair for $\Sigma = \{\symb{0},\symb{1}\}$ where the difference set $F= \llbracket -2,1\rrbracket$ is shown in boxes. Here $y = \sigma^{2}(x)$.
	
	\subsection{Basic properties of indistinguishable asymptotic pairs}

    We defined indistinguishable asymptotic pairs through patterns whose support is an arbitrary subset of $\ZZ$. 
    Next, we show that we can equivalently characterize indistinguishability using factors, first recalling the definitions.

     A pattern $w$ whose support is the set $\llbracket 0,n-1\rrbracket$ for some $n \in \NN$ is a \textbf{word}, and we write $w = w_0\dots w_{n-1}$. The length of $w$ is denoted by $|w|=n$.
    Let us denote the set of all finite words with symbols in $\Sigma$ by $\Sigma^* = \bigcup_{n \in \NN}\Sigma^{\llbracket 0,n-1\rrbracket}$.
        
    A word $w$ is a \textbf{factor} of $x$ if it appears in $x$.
    We write $\Lcal_n(x)$ for the set of all factors of $x$ of length $n$ and
    the \define{language of $x$} is the union $\Lcal(x)$ of the sets $\Lcal_n(x)$ for every $n \in \NN$.

		\begin{proposition}\label{prop:trivialite}
			An asymptotic pair $x,y \in \Sigma^{\ZZ}$ is indistinguishable if and only if for every $w \in \Sigma^*$ we have \[ \Delta_w(x,y)  = 0. \]
		\end{proposition}
		
		\begin{proof}
			One direction is obvious. Let us suppose that for every $w \in \Sigma^*$ we have $\Delta_w(x,y)  = 0$ and let $S \subseteq \ZZ$ be a support and $p \in \Sigma^S$. For every $m \in \ZZ$ we can define $p' \in \Sigma^{m+S}$ by $p'(m+s) = p(s)$ for every $s \in S$. For every sequence $x \in \Sigma^{\ZZ}$ we have that $x \in [p]$ if and only if $\sigma^{-m}(x)\in [p']$. Consequently we have that $\occ_p(x) = \occ_{p'}(x)+m$ and thus,
			
			\[ \Delta_{p'}(x,y) = \Delta_p(x,y) \mbox{ for every asymptotic pair } x,y \in \Sigma^{\ZZ}.  \]
			
			By the former argument, without loss of generality, we may assume that $S \subseteq \llbracket 0,n-1 \rrbracket$ for some large enough $n$.
		
			Notice that $[p]$ is the disjoint union of all $[w]$ where $w$ is a word of length $n$ such that $w|_S = p$. It follows that for any $z \in \Sigma^{\ZZ}$ we have $\indicator{[p]}(z) =  1$ if and only if there is a unique such $w$ such that $[w] \subseteq [p]$ and $\indicator{[w]}(z)=1$. Letting $F$ be the difference set of $x,y$ we obtain,
			
			\begin{align*}
			\Delta_p(x,y) & = \sum_{n \in F - S}\indicator{[p]}(\sigma^n(y))-\indicator{[p]}(\sigma^n(x)) \\
			& = \sum_{n \in F -  \llbracket 0,n-1 \rrbracket}\indicator{[p]}(\sigma^n(y))-\indicator{[p]}(\sigma^n(x)) \\
			& = \sum_{n \in F -  \llbracket 0,n-1 \rrbracket}  \sum_{\substack{w \in \Sigma^{ \llbracket 0,n-1 \rrbracket} \\ [w] \subseteq [p]} }\indicator{[w]}(\sigma^n(y))-\indicator{[w]}(\sigma^n(x)).
			\end{align*}
			Exchanging the order of the sums yields 
			
			\[ \Delta_p(x,y)  = \sum_{\substack{w \in \Sigma^{ \llbracket 0,n-1 \rrbracket} \\ [w] \subseteq [p]}} \Delta_{w}(x,y) = 0.  \]
			And thus $x,y$ is an indistinguishable asymptotic pair.
			\end{proof}
			
			Given a sequence $x \in \Sigma^{\ZZ}$, define its \define{reversal} $x^{R}$ by setting $x^{R}(n) = x(-n)$ for every $n \in \NN$. The next proposition states that indistinguishable asymptotic pairs are stable under actions of the affine group of $\ZZ$.
			
			\begin{proposition}\label{prop:shifted_SI-for-Z}
				Let $x,y\in \Sigma^{\ZZ}$ be an indistinguishable asymptotic pair. 
				\begin{enumerate}
					\item For every $n \in \ZZ$,  $(\sigma^n(x),\sigma^n(y))$ is an  indistinguishable asymptotic pair.
					\item $(x^R,y^R)$ is an indistinguishable asymptotic pair.
				\end{enumerate}
			\end{proposition}
			
			\begin{proof}
				Given a pattern $p$ with support $S$, denote by $-p$ the pattern with support $-S$ such that $-p(-s)=p(s)$. It is clear that if $x \in \Sigma^{\ZZ}$ and $n \in \NN$, then for every pattern $p$, \[ \occ_{p}(x) = \occ_p(\sigma^n(x))-n\quad \mbox{ and }\quad \occ_p(x) = -\occ_{-p}(x^R).   \]
				Note that if $x,y$ is an asymptotic pair with difference set $F$, then the difference set for $\sigma^n(x),\sigma^n(y)$ is $F-n$ and the difference set for $x^R,y^R$ is $-F$.
				
				From the relations on the occurrence sets, we obtain that for every pattern $p$ we have that \[ \Delta_p(x,y) = \Delta_p(\sigma^n(x),\sigma^n(y)) = \Delta_{-p}(x^R,y^R).  \]
				In particular, $x,y$ is indistinguishable if and only if $\sigma^n(x),\sigma^n(y)$ is indistinguishable if and only if $x^R,y^R$ is indistinguishable.
			\end{proof}
	
	Let us recall that a sequence $(x_m)_{m \in \NN}$ of sequences in $\Sigma^{\ZZ}$ converges to $\bar{x} \in \Sigma^{\ZZ}$ if for every $n \in \ZZ$ we have that $(x_m)_{n} = \bar{x}_{n}$ for all large enough $m \in \NN$. If $(x_m,y_m)_{m \in \NN}$ is a sequence of asymptotic pairs, it is natural to ask that both $(x_m)_{m \in \NN}$ and $(y_m)_{m \in \NN}$ converge to say that $(x_m,y_m)_{m \in \NN}$ converges. However, if we only asked for that there would be no guarantee that the limit is also an asymptotic pair. We shall consider a slightly stronger notion of convergence for asymptotic pairs which ensures that the limit is also an asymptotic pair.
	
	\begin{definition}
		We say that a sequence $(x_n,y_n)_{n \in \NN}$ of asymptotic pairs \textbf{converges} to a pair $(x,y)$ if $(x_n)_{n \in \NN}$ converges to $x$, $(y_n)_{n \in \NN}$ converges to $y$, and there exists a finite set $F \subseteq \ZZ$ so that $x_n|_{\ZZ \setminus F} = y_n|_{\ZZ \setminus F}$ for all large enough $n \in \NN$.
	\end{definition}
	
	This notion of convergence is also used in the theory of topological orbit equivalence of Cantor minimal systems. An interested reader can refer to~\cite{Put18} for further information. The advantage of this notion in our context is that it preserves indistinguishability.
	 
    \begin{proposition}\label{prop:limit-of-indist-is-indist}
		Let $(x_n,y_n)_{n \in \NN}$ be a sequence of asymptotic pairs in $\Sigma^{\ZZ}$ which converges to $(x,y)$. If for every $n \in \NN$ we have that $(x_n,y_n)$ is indistinguishable, then $(x,y)$ is indistinguishable.
	\end{proposition}
	
	\begin{proof}
		Let $p \in \Sigma^S$ be a pattern. As $(x_n,y_n)_{n \in \NN}$ converges to $(x,y)$, there exists a finite set $F \subseteq \ZZ$ and $N_1 \in \NN$ so that $x_n|_{\ZZ \setminus F} = y_n|_{\ZZ \setminus F}$ for every $n \geq N_1$. In particular we have that the difference sets of $(x,y)$ and $(x_n,y_n)$ for $n \geq N_1$ are contained in $F$. It suffices thus to show that \[ \# \{ \occ_p(x) \cap (F-S)\} = \#\{ \occ_p(y) \cap (F-S)\}.   \]
		As $(x_n)_{n \in \NN} $ converges to $x$ and $(y_n)_{n \in \NN} $ converges to $y$, there exists $N_2 \in \NN$ such that \[x_n|_{(F-S)+S} = x|_{(F-S)+S} \mbox{ and } y_n|_{(F-S)+S} = y|_{(F-S)+S} \mbox{ for all } n \geq N_2.\] This implies that $\occ_p(x) \cap (F-S) = \occ_p(x_n) \cap (F-S)$ and $\occ_p(y) \cap (F-S) = \occ_p(y_n) \cap (F-S)$ for every $n \geq N_2$.
		Taking $N = \max\{N_1,N_2\}$, as $(x_n,y_n)$ is indistinguishable with the difference set contained in $F$, it follows that for $n \geq N$ we have $\# \{ \occ_p(x_n) \cap (F-S)\} = \#\{ \occ_p(y_n) \cap (F-S)\}.$ Therefore we obtain $\# \{ \occ_p(x) \cap F-S\} = \#\{ \occ_p(y) \cap F-S\}.$ As this argument holds for every pattern $p$, we conclude that $(x,y)$ is indistinguishable.
	\end{proof}

    \subsection{Recurrence of indistinguishable asymptotic pairs}
	
    In this section, we study the recurrence of indistinguishable asymptotic pairs.
    To that end, we shall first show that if $x,y$ is an indistinguishable asymptotic pair, then every pattern with support $I\subseteq\ZZ$ which appears in $x$ must necessarily appear at some position $n \in \ZZ$ so that $n+I$ intersects the difference set of $x,y$. More precisely, let us say that the \define{occurrences of a word $w \in \Lcal(x)$ intersect $F$ in $x$} if $\occ_w(x) \cap (F-\llbracket 0,|w|-1\rrbracket) \neq \varnothing$. Equivalently, there exist $i \in F$ and $j \in \llbracket 0,|w|-1\rrbracket$ such that $\sigma^{i-j}(x) \in [w]$.
	
	\begin{lemma}\label{lem:words_appear_in_diff_set}
		Let $x,y \in \Sigma^{\ZZ}$ be a non-trivial indistinguishable asymptotic pair with the difference set $F$. The occurrences of every $w \in \Lcal(x)$ intersect $F$ in $x$.
	\end{lemma}
	
	\begin{proof}
		Without loss of generality, let us suppose that $F$ is contained in $\llbracket 0,k-1\rrbracket$, $x_0 \neq y_0$, and $x_{k-1} \neq y_{k-1}$. 
		Let us write $I = \llbracket 0,|w|-1\rrbracket$. If $(F-I) \cap \occ_w(x) = \varnothing$, then there is $u \in \occ_w(x)$ such that either (1) $u$ is the smallest value satisfying $u \geq k$ or (2) $u$ is the largest value satisfying $u \leq -|w|$. We shall deal with case the (1), the second case is analogous.
		
		As $\llbracket 0,k-1\rrbracket$ contains the difference set of $x,y$ we have that $x_{u}\dots,x_{u+|w|-1} = y_{u}\dots,y_{u+|w|-1}$. Let $w' = y_{k-1}\dots y_{u+|w|-1}$. As $x_{k-1}\neq y_{k-1}$, we obtain $\sigma^{k-1}(y) \in [w']$ but $\sigma^{k-1}(x) \notin [w']$. As $\Delta_{w'}(x,y)=0$, there must be some $j' \in F$ and $i' \in \llbracket 0,u+|w|-k\rrbracket$ such that $j'-i' \neq k-1$ and $\sigma^{j'-i'}(x) \in [w']$ and, since $w$ is a suffix of $w'$, $\sigma^{j'-i'+u -(k-1)}(x)\in [w]$.
		
		Let $u' = j'-i'+u-(k-1)$. On one hand, we have $j'-i' < k-1$ and thus $u' < u$.
		On the other hand, $u' \geq 0-(u+|w|-k)+u-(k-1) = -|w|+1$. 
		Since $u$ is the smallest value of $\occ_w(x)$ satisfying $u \geq k$, we obtain $u' \in F-I$ which is a contradiction.
	\end{proof}

	A sequence $x \in \Sigma^{\ZZ}$ is called \define{recurrent} if every $w \in \Lcal(x)$ occurs at least twice in $x$. It is quite easy to see that $x$ is recurrent if and only if $\occ_w(x)$ is in fact an infinite set for every $w \in \Lcal(x)$.
	
	We say that $x$ is \define{uniformly recurrent} if every $w \in \Lcal(x)$ appears with bounded gaps, that is, for every $w \in \Lcal(x)$ there exists an integer $g \geq 1$ such that for every $n \in \ZZ$ there is $0 \leq m \leq g$ such that $\sigma^{n+m}(x) \in [w]$.
	
	It is clear that if $(x,y)$ is an indistinguishable asymptotic pair, then $x$ is (uniformly) recurrent if and only if $y$ is (uniformly) recurrent.
	
	\begin{proposition}\label{prop:recurrence}
		Let $x,y \in \Sigma^{\ZZ}$ be an indistinguishable asymptotic pair. If $x$ is not recurrent, then $x$ and $y$ lie in the same orbit.
	\end{proposition}
	
	\begin{proof}
		If $x$ is not recurrent, there is a word $w \in \Lcal(x)$ such that $\#(\occ_w(x))=1$. Without loss of generality let us assume that $w$ occurs at the origin, that is, $\occ_w(x) = \{0\}$. As $\Delta_w(x,y)=0$, it follows that $\#(\occ_w(y))=1$ as well. Let $m$ be the only integer such that $\sigma^m(y) \in [w]$.
		
		Let $n\in \NN$ be larger than the length of $w$. Let $q_n = x|_{\llbracket -n,n \rrbracket}$. As $x \in [q_n]$ and $\Delta_{q_n}(x,y)=0$, there exists $k \in \ZZ$ so that $\sigma^{k}(y) \in [q_n]$. Furthermore, as $q_n|_{\llbracket 0,|w|-1\rrbracket} = w$, it follows that $\sigma^{k}(y) \in [w]$ and thus $k = m$. Therefore we obtain that $\sigma^{m}(y) \in [q_{n}]$ for every large enough $n$. As $\bigcap_{n \in \NN}[q_n] = \{x\}$ we deduce that $\sigma^{m}(y) = x$.
	\end{proof}
	
	Next we are going to show that recurrent indistinguishable asymptotic pairs are in fact uniformly recurrent. To that end, we recall the notions of return word~\cite{MR1489074,MR1829231,MR1808196} and 
	complete return word~\cite{MR2489283}.
	
	\begin{definition}\label{def:complete-return-word}
		A word $u \in \Sigma^+$ is a \define{complete return word} to $w \in \Sigma^+$
		in $x \in \Sigma^{\ZZ}$ if $u$ appears in $x$,
		$u=ws=pw$ for some nonempty words $p,s\in\Sigma^+$,
		and there are only two occurrences of $w$ in $u$, one as a prefix and
		one as a suffix. 
		The word $p$ is called a \define{return word} to $w$ in $x$.
	\end{definition}
	
	Note that the two occurrences of $w$ in a complete return word $u=ws=pw$ to $w$ may overlap.
	Denote the set of all complete return words to $w$ in $x$ by $\texttt{CRW}_{w}(x)$ and the set of all return words to $w$ in $x$ by $\texttt{RW}_w(x)$. The following fact is elementary.
	
	\begin{lemma}\label{lem:CWR_finite_iff_unifrec}
		Let $x \in \Sigma^{\ZZ}$. The following are equivalent.
		\begin{enumerate}
			\item $x$ is uniformly recurrent.
			\item $x$ is recurrent and for every $w \in \Lcal(x)$ we have that $\texttt{CRW}_{w}(x)$ is finite.
			\item $x$ is recurrent and for every $w \in \Lcal(x)$ we have that $\texttt{RW}_{w}(x)$ is finite.
		\end{enumerate}
	\end{lemma}
	
	\begin{lemma}\label{lem:rec_implies_urec}
		Let $x,y\in \Sigma^{\ZZ}$ be a non-trivial indistinguishable asymptotic pair. If $x$ is recurrent, then $x$ is uniformly recurrent.
	\end{lemma}
	
	\begin{proof}
		Without loss of generality, suppose the difference set of $(x,y)$ is contained in $F = \llbracket 0, k-1\rrbracket$. Suppose $x$ is not uniformly recurrent. By~\Cref{lem:CWR_finite_iff_unifrec} there is a word $w \in \Lcal(x)$ such that $\texttt{CRW}_{x}(w)$ is infinite. As $\texttt{CRW}_{x}(w)$ is infinite, there exist distinct $v_1,v_2,v_3 \in \texttt{CRW}_{x}(w)$ such that $\min\{|v_1|,|v_2|\}> k+2|w|$ and $|v_3|> k+2\max\{|v_1|,|v_2|\}$. By~\Cref{lem:words_appear_in_diff_set} the words $v_1,v_2$ and $v_3$ must occur in $x$ at positions such that their support intersects $F = \llbracket 0, k-1\rrbracket$. 
		
		As $\min\{|v_1|,|v_2|\}> k+2|w|$, exactly one of the two occurrences of $w$ in $v_1$ must be completely contained in the support $ L_1 = \llbracket -|v_1|+1,-1\rrbracket$ or the support $R_1 = \llbracket k,k+|v_1|-1\rrbracket$. Similarly, exactly one occurrence of $w$ in $v_2$ appears in  $L_2 = \llbracket-|v_2|+1,-1\rrbracket$ or $R_2=\llbracket k,k+|v_2|-1\rrbracket$. As $v_1,v_2$ are distinct complete return words, if an occurrence of $w$ coming from $v_1$ appears in $L_1$, then another coming from $v_2$ appears in $R_2$. Analogously, if there is an occurrence of $w$ coming from $v_1$ in $R_1$, then an occurrence of $w$ coming from $v_2$ appears in $L_2$.
		
		In consequence with the reasoning above, the word $w$ appears completely contained both in the interval $\llbracket-\max\{|v_1|,|v_2|\}+1,-1\rrbracket$ and in the interval $\llbracket k,k+\max\{|v_1|,|v_2|\}-1\rrbracket$. As $v_3$ is also a complete return word which appears intersecting $F$ and $|v_3|> k+2\max\{|v_1|,|v_2|\}$, this means that there are no copies of $w$ completely contained in either $\llbracket-\max\{|v_1|,|v_2|\}+1,-1\rrbracket$ or $\llbracket k,k+\max\{|v_1|,|v_2|\}-1\rrbracket$, contradicting the above statement. \end{proof}
	
	Gathering~\Cref{prop:recurrence} and~\Cref{lem:rec_implies_urec} we obtain the following beautiful dichotomy.
	
	\begin{corollary}
		Let $x,y\in \Sigma^{\ZZ}$ be a non-trivial asymptotic indistinguishable pair. Then exactly one of the following statements holds \begin{enumerate}
			\item $x = \sigma^n(y)$ for some nonzero $n \in \ZZ$,
			\item $x$ and $y$ are uniformly recurrent.
		\end{enumerate}
	\end{corollary}
	
	\begin{proof}
		If $x$ is not recurrent, then by~\Cref{prop:recurrence} we obtain that $x = \sigma^n(y)$ for some $n \in \ZZ\setminus \{0\}$. If $x$ is recurrent, then $y$ is also recurrent. Applying~\Cref{lem:rec_implies_urec} we obtain that $x$ and $y$ are uniformly recurrent.
		
		Let us assume that both conditions happen at the same time. As $x = \sigma^n(y)$ for some nonzero $n \in \ZZ$ and $x,y$ are asymptotic, we obtain that $x$ is eventually periodic. Furthermore, as $x$ is uniformly recurrent, we obtain that $x$ is a periodic sequence. Hence the only possibility to have a finite difference set is having $x=y$, which contradicts the non-triviality assumption.
	\end{proof}
	
	\section{Lower and upper characteristic Sturmian sequences on $\ZZ$}\label{sec:sturmian}

    The purpose of this section is to
    prove~\Cref{thm:sturmian_characterization}, that is,
    that recurrent indistinguishable asymptotic pairs $x,y\in\{\symb{0},\symb{1}\}^\ZZ$ 
    whose difference set is of size 2
    consist of lower and upper characteristic Sturmian sequences and vice versa.
    The proof of the first implication is based on the description of Sturmian sequences
    as lower and upper mechanical words, a terminology introduced by
    Morse and Hedlund \cite{MR0000745}. Notice that here the word ``mechanical
    word'' is used to refer to a biinfinite sequence in our context.
    The proof of the reciprocal is based on the description of Sturmian sequences
    by their factor complexity \cite{MR0322838}.
    Schematically, the proofs in this section 
    are done as follows:
    \begin{center}
        \includegraphics{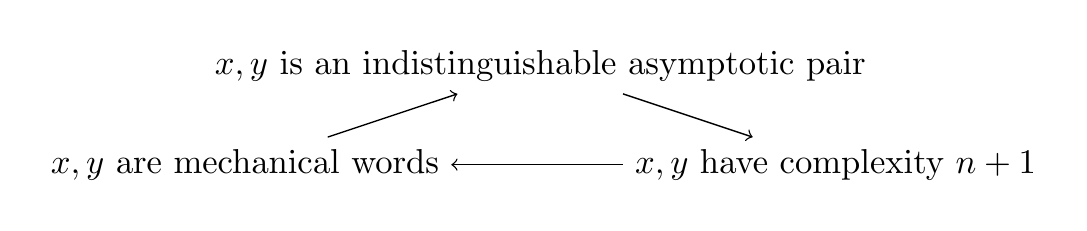}
    \end{center}
    The fact that
    recurrent sequences with
    complexity $n+1$
    are
    mechanical words of irrational slope
    implies that all the aforementioned are equivalent.

    \subsection{Mechanical words}
    Given two real numbers $\alpha$ and $\rho$ with $0\leq\alpha< 1$, we define two sequences
	\[
	s_{\alpha,\rho}:\ZZ\to \{\symb{0},\symb{1}\},\quad
	s'_{\alpha,\rho}:\ZZ\to \{\symb{0},\symb{1}\}
	\]
	by
	\begin{align*}
	s_{\alpha,\rho}(n) = \lfloor\alpha(n+1)+\rho\rfloor-\lfloor\alpha n+\rho\rfloor,\\
	s'_{\alpha,\rho}(n) = \lceil\alpha(n+1)+\rho\rceil-\lceil\alpha n+\rho\rceil.
	\end{align*}
	The sequence $s_{\alpha,\rho}$ is the \define{lower mechanical word} and
	$s'_{\alpha,\rho}$ is the \define{upper mechanical word} with \define{slope} $\alpha$
	and \define{intercept} $\rho$, see Chapter 2 of~\cite{MR1905123}.
	It is clear that if $\rho-\rho'$ is an integer, then $s_{\alpha,\rho}=s_{\alpha,\rho'}$ and	$s'_{\alpha,\rho}=s'_{\alpha,\rho'}$. Thus we may always assume $0\leq\rho<1$.

	The mechanical words $s_{\alpha,\rho},s'_{\alpha,\rho}$ are in fact codings of trajectories of irrational circle rotations, namely, consider the isometry $R_\alpha\colon \RR/\ZZ \to \RR/\ZZ$, where $R_\alpha(\rho) = (\rho + \alpha) \bmod 1$ for every $\rho \in \RR/\ZZ$.
	Consider the partition $\mathcal{P} = \{ I_{\symb{0}}, I_{\symb{1}}\}$ of $\RR/\ZZ$ given by $I_{\symb{0}} = [0,1-\alpha)$ and $I_{\symb{1}} = [1-\alpha,1)$.
	For $\rho \in \RR/\ZZ$, define
	\[
	\nu(\rho) = i \quad \text{ if } \rho \in I_i \quad \text{ for } i \in \{\symb{0},\symb{1}\}.
	\]
	We obtain
	\[
	s_{\alpha,\rho}(n) = \nu(R_{\alpha}^n(\rho)) \quad \mbox{ for every } n \in \ZZ,
	\]
	i.e., $s_{\alpha,\rho}$ is the coding of the trajectory of $\rho$ with respect to the partition $\mathcal{P}$, see Section 2.2.2 of~\cite{MR1905123}. Similarly, $s_{\alpha,\rho}$ is the coding of the trajectory of $\rho$ with respect to the partition  $\mathcal{P}' = \{ I'_{\symb{0}}, I'_{\symb{1}}\}$ of $\RR/\ZZ$ given by $I'_{\symb{0}} = (0,1-\alpha]$ and $I'_{\symb{1}} = (1-\alpha,1]$.

	Since $1+\lfloor\alpha n+\rho\rfloor=\lceil\alpha n+\rho\rceil$ whenever $\alpha n+\rho$
	is not an integer, one has 
	$s_{\alpha,\rho}=s'_{\alpha,\rho}$ except when $\alpha n+\rho$ is an
	integer for some $n\in\ZZ$. As $\alpha \notin \QQ$, this can happen for at most one $n \in \ZZ$, in this case,
	\begin{align*}
	&s_{\alpha,\rho}(n-1) = \symb{1},
	&&s'_{\alpha,\rho}(n-1) = \symb{0}, \\
	&s_{\alpha,\rho}(n) = \symb{0},
	&&s'_{\alpha,\rho}(n) = \symb{1},
	\end{align*}
	and elsewhere
	\[
	s_{\alpha,\rho}(k) = s'_{\alpha,\rho}(k) \qquad\text{ whenever }\quad k\notin\{n-1,n\}.
	\]

    In what follows, we say that
    the sequence $c_\alpha=s_{\alpha,0}$ is the \define{lower
    characteristic Sturmian sequence of slope $\alpha$} and the sequence $c'_\alpha=s'_{\alpha,0}$
    is the \define{upper characteristic Sturmian sequence of slope $\alpha$}. Notice that
    $c_{\alpha}(n)=c'_{\alpha}(n)$ if and only if $n\in\ZZ\setminus\{-1,0\}$. 

    \begin{remark}\label{rem:two-sided-characteristic}
    For one-sided sequences, the characteristic Sturmian sequence of
    slope $\alpha$ is usually the one having two
    distinct extensions to the left, see \cite{MR1902372}, \cite{MR1466181} or
        \cite[\S 9]{MR1997038}.
    Here, we consider biinfinite Sturmian sequences as in 
        \cite[\S 6.2]{MR1970391}, and we believe it is more natural to
    define $s_{\alpha,0}$ and $s'_{\alpha,0}$ with intercept $\rho=0$ as the
        lower and upper ``characteristic" ones.
    \end{remark}

    \subsection{Pairs of characteristic Sturmian sequences are indistinguishable}

	\begin{proposition}\label{prop:Sturmian_are_SI}
		The lower and upper characteristic Sturmian sequences
        $(c_{\alpha},c'_{\alpha})$ 
        form a non-trivial indistinguishable asymptotic pair for
        every irrational $\alpha\in[0,1]\setminus\QQ$. 
	\end{proposition}

	\begin{proof}
		From the above discussion, it follows that $c_{\alpha}$
		and $c'_{\alpha}$ are asymptotic with the difference set $\left\{-1, 0\right\}$.
        The pair is non-trivial since the difference set is nonempty.
		Note that in this case,
		\begin{align*}
		&c_{\alpha}(-1) = \symb{1},
		&&c'_{\alpha}(-1) = \symb{0}, \\
		&c_{\alpha}(0) = \symb{0},
		&&c'_{\alpha}(0) = \symb{1}.
		\end{align*}
		
		Let $m \in \NN$ and $w \in \{\symb{0},\symb{1}\}^m$.
		We shall show that 
\begin{equation} \label{eq:St_are_SI_eq}
    \sum_{i=0}^m \indicator{[w]}\sigma^{-i}(c_{\alpha}) = 
    \sum_{i=0}^m \indicator{[w]}\sigma^{-i}(c'_{\alpha}).
\end{equation}

		Note that this sum above has $m+1$ indexes. As $c_{\alpha},c'_{\alpha}$ are Sturmian of angle $\alpha$, we have that $\Lcal_m(c_{\alpha}) = \Lcal_m(c'_{\alpha})$ is of size $m+1$.
		Together with showing that for each word $w \in \Lcal_m(c_{\alpha})$ there exist 
        $i,i' \in \llbracket 0,m\rrbracket$ 
        such that  $\sigma^{-i}(c_{\alpha}) \in [w]$ and $\sigma^{-i'}(c'_{\alpha}) \in [w]$, it implies that such $i$ and $i'$ are unique, which implies \eqref{eq:St_are_SI_eq}.
		
		Let us consider the refinement $\mathcal{P}^m = \bigvee_{j \in \llbracket 0,m-1\rrbracket} R_{\alpha}^{-j}(\mathcal{P})$. That is, the partition obtained by intersecting the semiclosed intervals of each shifted partition $R_{\alpha}^{-j}(\mathcal{P})$ between themselves. By definition of the coding $\nu$, for each $w \in \Lcal_m(c_{\alpha})$ there is $I \in \mathcal{P}^m$ such that for every $x \in I$ we have \[\nu(x)\nu(R_{\alpha}(x))\cdots\nu(R^{m-1}_{\alpha}(x)) = w. \]
		
		In consequence, there are $m+1$ semiclosed intervals in $\mathcal{P}^m$ representing each word in $\Lcal_m(c_{\alpha})$. Note that the set of (closed) endpoints of the semiclosed intervals in $\mathcal{P}^m$ is given by the collection:
		\[
		\left\{ 0, -\alpha \bmod{1}, -2\alpha \bmod{1},\dots,-m\alpha \bmod{1} \right\}.
		\]
		
		As for $i \in \llbracket 0,m\rrbracket$ we have $\sigma^{-i}(c_{\alpha}) = s_{\alpha,-i\alpha }$, we obtain that each of these shifts $\sigma^{-i}(c_{\alpha})$ begins in one of the above endpoints.
        This proves that there exists $i \in \llbracket 0,m\rrbracket$ such that  $\sigma^{-i}(c_{\alpha}) \in [w]$.
		
		The situation for the upper characteristic word $c'_{\alpha}$ is analogous with the following distinction: all intervals are left-open right-closed, that is, the initial partition is $I_{\symb{0}}' = (0,1-\alpha]$ and $I_\symb{1}' = (1-\alpha,1]$.
		The analogous partition $(\mathcal{P}')^m$ has the same set of endpoints as $\mathcal{P}^m$ and thus the same conclusion follows.

		Equality \eqref{eq:St_are_SI_eq} implies $\Delta_w(x,y) = 0$ for all $w \in \Lcal_m(c_{\alpha})$ and all $m$.
		By \Cref{prop:trivialite}, the lower and upper characteristic words form an indistinguishable pair.		
	\end{proof}
	
	\begin{remark}
		If $\{\alpha n+\rho : n\in\ZZ\}\cap\ZZ=\varnothing$, then
		$s_{\alpha,\rho} =s'_{\alpha,\rho}$ and then $s_{\alpha,\rho}$ and $s'_{\alpha,\rho}$ form a trivial indistinguishable asymptotic pair.
		Otherwise, if there exists $n \in \ZZ$ such that $\alpha n+\rho \in \ZZ$, then $c_{\alpha} = \sigma^{n}(s_{\alpha,\rho})$ and $c'_{\alpha} = \sigma^{n}(s'_{\alpha,\rho})$.
		By~\Cref{prop:Sturmian_are_SI,prop:shifted_SI-for-Z} it follows that $s_{\alpha,\rho}$ and $s'_{\alpha,\rho}$ form a non-trivial indistinguishable asymptotic pair.
	\end{remark}

\subsection{Recurrent indistinguishable asymptotic pairs are Sturmian}

The goal of this subsection is to prove the reciprocal
about recurrent indistinguishable asymptotic pairs with difference set $F =
\{-1,0\}$. In order to do that, we shall first show that their factor complexity,
which counts the number of words of each length in their language,
coincides with that of a Sturmian sequence.

    The \define{factor complexity} of a sequence $x \in \{\symb{0},\symb{1}\}^{\ZZ}$
    is the mapping $n \mapsto \# \Lcal_n(x)$.
    Let us recall that a sequence $x \in \{\symb{0},\symb{1}\}^{\ZZ}$ is \define{Sturmian} if its factor complexity $\#\Lcal_n(x) = n+1$ for every $n \in \NN$ and it is not eventually periodic \cite[Def.~6.2.4, Prop.~6.2.5]{MR1970391}. 
    Moreover, a sequence is Sturmian if and only if it is a lower or upper
    mechanical word for some irrational slope $\alpha$
    \cite{MR0000745,MR0322838}, see also \cite[Theorem 2.1.13]{MR1905123}.

    The study of factor complexity is closely related to special factors,
    a notion which is used in the next proof to provide a lower bound.
    A word $w\in\Lcal_n(x)$ is called \define{right special} (\define{left special} resp.)
    in $x$ if there exists at least two distinct letters
    $a,b\in\Sigma$ such that $wa,wb\in\Lcal_{n+1}(x)$ (such that
    $aw,bw\in\Lcal_{n+1}(x)$ resp.), see \cite{MR2759107}.
	
	A consequence of \Cref{lem:words_appear_in_diff_set} is that the factor complexity of indistinguishable asymptotic pairs can be bounded above by the size of the smallest interval that contains their difference set.

	\begin{proposition}\label{prop:complexity-bound}
		Let $x,y \in \Sigma^{\ZZ}$ be a
		non-trivial indistinguishable asymptotic pair whose difference set $F$ is contained in an interval $I$.
		We have that for every $n \geq 1$ \[n+1 \leq \#\Lcal_n(x) \leq n+\#(I)-1.\]
	\end{proposition}
	
	\begin{proof}
		By~\Cref{lem:words_appear_in_diff_set}, the occurrences of every $w \in \Lcal_n(x)$ intersect $F$ in $x$.
		In other words, for each $w \in \Lcal_n(x)$ there exists $u \in F - \llbracket 0,n-1\rrbracket$ such that $\sigma^u(x) \in [w]$.
		Without loss of generality, by shifting $x$ and $y$ we may assume that $F \subseteq \llbracket 0,k-1\rrbracket$ (hence $I = \llbracket 0,k-1\rrbracket$), so there exists a surjective function from $\llbracket -n+1, k-1 \rrbracket$ to $\Lcal_n(x)$. In particular, $\#\Lcal_n(x)\leq n+k-1=n+\#(I)-1$.
		
		In order to obtain the lower bound, notice that as $x \neq y$, we have $F \neq \varnothing$ and thus we can define $m=\max\{i\in\ZZ\mid x_i\neq y_i\}$.
		For all $n\geq0$, we have
		that the word \[w=x_{m+1}\dots x_{m+n} = y_{m+1}\dots y_{m+n}\] can be left-extended to a word of length $n+1$ in $\Lcal_{n+1}(x)$ in two different ways, namely \[ w' = x_mw, \quad w'' = y_mw. \]
        Thus, $w$ is a left special factor in $x$. Since every factor in $\Lcal_n(x)$ can be extended to the left by one symbol to get a word in $\Lcal_{n+1}(x)$
		and for every $n$ there exists a 
        left special factor of length $n$ in $x$,
		it implies that $\#\Lcal_{n+1}(x)-\#\Lcal_n(x)\geq1$ for every $n\geq0$.
		Since $x\neq y$, then $\#\Lcal_1(x)=\#\Sigma\geq 2$ and we conclude $\#\Lcal_n(x) \geq n+1$.
	\end{proof}

As a consequence, when the difference set of $x$ and $y$ is of size 2, the
factor complexity must be $n+1$.

	\begin{corollary}\label{cor:complexity-is-n+1}
        If $x,y \in \Sigma^{\ZZ}$ is a
        non-trivial indistinguishable asymptotic pair with difference set $F=\{-1,0\}$,
        then $\#\Lcal_n(x) = n+1$.
	\end{corollary}
	
    \begin{proof}
        By~\Cref{prop:complexity-bound} we deduce that 
        $n + 1 \leq \#\Lcal_n(x) = \#\Lcal_n(y) \leq n+ \# F -1=n+1$ for every $n \in \NN$ and thus 
        $\#\Lcal_n(x) = \#\Lcal_n(y) = n+1$ for every $n \in \NN$. 
    \end{proof}

It is known that for each Sturmian sequence and each nonnegative
integer $n$, some factor of length $2n$ of the sequence contains the $n + 1$
factors of length $n$ of the sequence, see for instance 
\cite[Corollary 5.2]{MR2197281}.
It turns out that the central factors of $x$ and $y$ of length $2n$ provide two
such words.

    \begin{corollary}\label{cor:two-compact-representation-of-the-language}
        If $x,y \in \Sigma^{\ZZ}$ is a
        non-trivial indistinguishable asymptotic pair with difference set $F=\{-1,0\}$,
        then each of the words
        $x_{-n}x_{-n+1}\cdots x_{n-1}$
        and
        $y_{-n}y_{-n+1}\cdots y_{n-1}$
        contain exactly one occurrence of each word in $\Lcal_n(x)$.
	\end{corollary}

    \begin{proof}
        From \Cref{cor:complexity-is-n+1},
        $\#\Lcal_n(x) = \#\Lcal_n(y) = n+1$ for every $n \in \NN$. 
        From \Cref{lem:words_appear_in_diff_set},
        both $x_{-n}x_{-n+1}\cdots x_{n-1}$
        and $y_{-n}y_{-n+1}\cdots y_{n-1}$ contain an occurrence of every factor
        of $\Lcal_n(x)$. All of the occurrences must be distinct or otherwise
        $\#\Lcal_n(x) < n+1$.
    \end{proof}

    For example,
    the following two words of length 26
    contains the same 14 factors of length 13:
    \begin{align*}
        &\symb{1010010100101\zeropoint0010010100101}\\
        &\symb{1010010100100\zeropoint1010010100101}
    \end{align*}

    It is well-known that
    one-sided sequences of complexity $n+1$
    are not eventually periodic, see
    \cite[Th.~6.1.8]{MR1970391}
    and
    \cite[Th.~2.1.13]{MR1905123}.
    This is no longer true for biinfinite sequences of complexity $n+1$,
    e.g., consider 
    $\leftidx^{\infty}\symb{0}\zeropoint\symb{1}^{\infty}$
    or $\leftidx^{\infty}\symb{0}\zeropoint\symb{1}\symb{0}^{\infty}$.
    A way to exclude eventually periodic sequences of complexity $n+1$ in the biinfinite setting
    is to consider recurrent sequences.
    For completeness, we provide a proof of the following result which can be
    considered as folklore even if not mentioned in \cite[\S 6.2]{MR1970391}.

\begin{proposition}\label{prop:biinfinite-sturmian} $x\in\{\symb{0},\symb{1}\}^\ZZ$ is Sturmian
    if and only if
    $x$ is recurrent and $\#\Lcal_n(x)=n+1$.
\end{proposition}

\begin{proof}
Sturmian sequences are recurrent, see \cite[Exercise 6.2.10]{MR1970391}.

Assume now that $x$ is recurrent and $\#\Lcal_n(x)=n+1$.
For contradiction, assume that $x$ is eventually periodic. 
Let $x = vp^\infty$ where $p$ is the shortest such word (and $v$ is a one-sided left infinite word). %
The choice of $p$ implies that the set of factors of length $|p|$ of the word $p^\infty$ has exactly $|p|$ elements.
As $\Lcal_{|p|}(x) = |p| + 1$, there is a word $u \in \Lcal(x)$ which does not occur in $p^\infty$, and there is a last occurrence of $u$ in $x$.
The last occurrence of $u$ in $x$ is followed by an arbitrarily long factor $s$ with no occurrence of $u$, and, as $x$ is recurrent, the factor $us$ has infinitely many occurrences in $v$.
Therefore, $u$ has unbounded gaps between its occurrences.
    In particular, there are at least three distinct complete return words
    $\{r_1,r_2,r_3\}$ to $u$ in $\Lcal(x)$. More precisely, 
    for each $a\in\{1,2,3\}$,
    $r_a$ contains exactly two occurrences of $u$, one as a prefix and one as a suffix.
    Let $p_{ab}$ denote the longest common prefix of $r_a$ and $r_b$.
    Up to some permutation of the complete return words,
    we have 
    \[
        |p_{12}| = 
        |p_{13}| <
        |p_{23}|.
    \]
    The word $p_{12}$ is a right special factor in $\Lcal(x)$, that is, there exist
    two distinct letters $a,b\in\{\symb{0},\symb{1}\}$ such that $p_{12}a,p_{12}b\in\Lcal(x)$.
    The suffix $s$ of $p_{23}$ of length $|p_{12}|$ is also a right special factor in $\Lcal(x)$.
    Moreover, $p_{12}\neq s$ since $u$ is a prefix of $p_{12}$ and $s$ contains
    no occurrence of $u$. This implies that
    $\#\Lcal_{N+1}(x)-\#\Lcal_N(x)\geq2$ for some integer $N\in\NN$, a
    contradiction.

The case with $w$ being eventually periodic to the left, i.e., $w = {}^\infty pv$, is analogous.
We conclude that $x$ is not eventually periodic, and thus it is Sturmian by definition.
\end{proof}

    These complexity bounds are the main tools to provide the characterization
    of pairs consisting of lower and upper characteristic Sturmian sequences by recurrent
    indistinguishable asymptotic pairs. 

    \begin{proof}[Proof of~\Cref{thm:sturmian_characterization}]
        By~\Cref{prop:Sturmian_are_SI}, the lower characteristic word $c_{\alpha}$
		and the upper characteristic word $c'_{\alpha}$ 
        form an indistinguishable asymptotic pair for
        every irrational $\alpha$ with $F=\{-1,0\}$ as their difference set where
        $x_{-1}x_{0}=\symb{10}$
        and $y_{-1}y_{0}=\symb{01}$ for $x={c}_{\alpha}$ and $y={c}'_{\alpha}$.
            If $\alpha\in[0,1]\setminus\QQ$, then $x$ and $y$ are recurrent.

		Conversely, from \Cref{cor:complexity-is-n+1}
        we have $\#\Lcal_n(x)=n+1$.
        From \Cref{prop:biinfinite-sturmian},
        recurrent sequences on $\ZZ$ of complexity $n+1$ are Sturmian sequences.
        We conclude that $x$ and $y$ are Sturmian sequences with the same language
        associated to some irrational slope $\alpha\in[0,1]\setminus\QQ$.
        Since $x|_{\ZZ\setminus F}=y|_{\ZZ\setminus F}$ with
        $F=\{-1,0\}$ such that $x_{-1}x_{0}=\symb{10}$ and $y_{-1}y_{0}=\symb{01}$
        we conclude that $x={c}_{\alpha}$ and $y={c}'_{\alpha}$ 
        are respectively the lower and upper characteristic Sturmian sequences with slope $\alpha$.
	\end{proof}
	
    \begin{remark}\label{rem:toeplitz}
		One might wonder if it is possible to prove~\Cref{thm:sturmian_characterization} with weaker assumptions, for instance, by asking just that $\Lcal(x)=\Lcal(y)$ instead of indistinguishability. This particular condition does not suffice, even if we further ask that the sequences are uniformly recurrent. Indeed, let $z \in \{\symb{0},\symb{1}\}^{\ZZ \setminus \{0\}}$ be defined by $z(n) = k \bmod 2$ whenever $k \geq 1$ and $n = 2^{k-1} \bmod 2^k$. Notice that $z$ is well defined for every nonzero integer and looks as follows 
		\[ z = \dots \symb{1}\symb{0}\symb{1}\symb{1}\symb{1}\symb{0}\symb{1}\symb{0}\symb{1}\symb{0}\symb{1}\symb{1}\symb{1}\symb{0}\symb{1}\zeropoint \symb{?} \symb{1}\symb{0}\symb{1}\symb{1}\symb{1}\symb{0}\symb{1}\symb{0}\symb{1}\symb{0}\symb{1}\symb{1}\symb{1}\symb{0}\symb{1}\dots    \]
		Let $x,y \in \{\symb{0},\symb{1}\}^{\ZZ}$ be the asymptotic pair defined by \[ x(n) = \begin{cases}
		\symb{0} & \mbox{ if } n = 0\\
		z(n) & \mbox{ otherwise }\\
		\end{cases} \quad \mbox{ and }\quad y(n) = \begin{cases}
		\symb{1} & \mbox{ if } n = 0\\
		z(n) & \mbox{ otherwise }\\
		\end{cases}.  \]
		
		The sequences $x$ and $y$ are limits of \define{Toeplitz sequences} which were defined in~\cite{Jacobs196901sequencesOT}. They are uniformly recurrent (see for instance \cite[Section 4.6]{Ku_ToSyDy}), have the same language and are not Sturmian. Furthermore, if one wishes to construct an example with difference set $\{-1,0\}$, one can consider the Thue-Morse substitution $\varphi \colon \{\symb{0},\symb{1}\} \to \{\symb{0},\symb{1}\}^*$ given by $\varphi(\symb{0}) = \symb{0}\symb{1}$ and $\varphi(\symb{1})=\symb{1}\symb{0}$ and consider $x'=\sigma(\varphi(x))$ and $y'=\sigma(\varphi(y))$. Then $x',y'$ are uniformly recurrent, form an asymptotic pair with the same language, and have difference set $\{-1,0\}$. A direct inspection of their language shows that they are not Sturmian.

	\end{remark}

	\section{Limits of Sturmian sequences toward rational slopes}\label{sec:limit-sturmian}

    In this section, we describe
    limits of Sturmian sequences
    toward a rational slope from above or from below
    and we show that they also constitute
    indistinguishable asymptotic pairs in $\ZZ$.
    Such words were already considered for instance 
    in \cite{MR2242616} and
    in \cite{MR2387854} (see condition $B_4$ and Figure 2).
	We prove
    \Cref{thm:sturmian_characterization_non_recurrent}
    about non-recurrent sequences 
    which,
    together with \Cref{thm:sturmian_characterization},
    implies \Cref{thm:sturmian_characterization_Z1}.

    \subsection{Christoffel words}

Christoffel words have many equivalent definitions, see \cite{MR1453849,MR2681739}
and the books \cite{MR2464862,MR3887697}.
Let $p,q\in\ZZ$ be coprime integers
such that $p/q\in\QQ_{\geq0}\cup\{\infty\}$ where the limit cases are written
as $0=0/1$ and $\infty=1/0$.
The \define{lower Christoffel word of slope} $p/q$ is
the factor of length $p+q$ of the lower mechanical word of slope
$\alpha=p/(p+q)$ and intercept $\rho=0$ starting at index 0:
\[
    c_{\alpha}(0)
    c_{\alpha}(1)
    \cdots
    c_{\alpha}(p+q-1).
\]
Similarly,
the \define{upper Christoffel word of slope} $p/q$ is
the factor of length $p+q$ of the upper mechanical word of slope
$\alpha=p/(p+q)$ and intercept $\rho=0$ starting at index 0:
\[
    c'_{\alpha}(0)
    c'_{\alpha}(1)
    \cdots
    c'_{\alpha}(p+q-1).
\]
If $p=0$ and $q=1$, then the lower and upper Christoffel word of slope $p/q=0$ is $\symb{0}$.
If $p=1$ and $q=0$, then the lower and upper Christoffel word of slope $p/q=\infty$ is $\symb{1}$.

    \subsection{Limits of Sturmian sequences}

    The lower and upper characteristic Sturmian sequences are related to each other
    as one is the shifted reversal of the other.
    More precisely we have the following elementary result based on the symmetry of floor
    and ceiling functions.

    \begin{lemma}\label{lem:symmetry-lower-to-upper}
        The lower characteristic Sturmian sequence is the shifted reversal of the upper
        characteristic Sturmian sequence in the sense that
        $c_\alpha(n)=c'_\alpha(-n-1)$ for every $n\in\ZZ$.
        Moreover
        \[
        c_\alpha(n)=c_\alpha(-n-1) 
            \qquad
            \text{ and }
            \qquad
        c'_\alpha(n)=c'_\alpha(-n-1) 
        \]
        for every $n\in\ZZ\setminus\{-1,0\}$.
    \end{lemma}

    \begin{proof}
        For all $x\in\RR$, we have $\lfloor x\rfloor = -\lceil -x \rceil$.
        Thus
        $
        c_{\alpha}(n) 
            = \lfloor\alpha(n+1)\rfloor-\lfloor\alpha n\rfloor
            = \lceil\alpha(-n)\rceil-\lceil\alpha (-n-1)\rceil
            = c'_{\alpha}(-n-1)$.
    Also, if $n\in\ZZ\setminus\{-1,0\}$, then
        $ c_\alpha(n) =c'_{\alpha}(-n-1) =c_\alpha(-n-1)$.
    The same holds for $c'_\alpha$.
    \end{proof}

Limits of lower or upper characteristic Sturmian sequences toward rational slopes
are eventually periodic sequences of complexity $n+1$ which 
can be expressed in terms of Christoffel words.

\begin{lemma}\label{lem:limit-of-sturmian-as-christoffel}
    Let $p,q\in\ZZ_{\geq0}$ be coprime integers.
    Limits of lower or upper characteristic Sturmian sequences as their slope tends towards $p/(p+q)$ are of one of the following forms depending on the value of $p$ and $q$.
    If $p\neq0$ and $q\neq0$, then
    \begin{align*}
        \lim_{\alpha\to\frac{p}{p+q}^+} c_\alpha &= \leftidx{^\infty}(1m0)(1m1)\zeropoint(0m1)(0m1)^\infty, \\
        \lim_{\alpha\to\frac{p}{p+q}^+}c'_\alpha &= \leftidx{^\infty}(1m0)(1m0)\zeropoint(1m1)(0m1)^\infty, \\
        \lim_{\alpha\to\frac{p}{p+q}^-} c_\alpha &= \leftidx{^\infty}(0m1)(0m1)\zeropoint(0m0)(1m0)^\infty, \\
        \lim_{\alpha\to\frac{p}{p+q}^-}c'_\alpha &= \leftidx{^\infty}(0m1)(0m0)\zeropoint(1m0)(1m0)^\infty,
    \end{align*}
    where $0m1$ and $1m0$ are respectively the lower and upper Christoffel word of slope $p/q$
    with $m\in\{0,1\}^*$.
When $p=0$ and $q=1$ and the limit is done from above, then
    \begin{align*}
        \lim_{\alpha\to\frac{p}{p+q}^+} c_\alpha 
        &= \lim_{\alpha\to 0^+} c_\alpha = \leftidx{^\infty}01\zeropoint00^\infty,\\
        \lim_{\alpha\to\frac{p}{p+q}^+}c'_\alpha 
        &= \lim_{\alpha\to 0^+}c'_\alpha = \leftidx{^\infty}00\zeropoint10^\infty.
    \end{align*}
When $p=1$ and $q=0$ and the limit is done from below, then
    \begin{align*}
        \lim_{\alpha\to\frac{p}{p+q}^-} c_\alpha 
        &= \lim_{\alpha\to 1^-} c_\alpha = \leftidx{^\infty}11\zeropoint01^\infty, \\
        \lim_{\alpha\to\frac{p}{p+q}^-}c'_\alpha 
        &= \lim_{\alpha\to 1^-}c'_\alpha = \leftidx{^\infty}10\zeropoint11^\infty.
    \end{align*}
\end{lemma}

\begin{proof}
    Recall that $c_{\alpha}(n) = \lfloor\alpha(n+1)\rfloor-\lfloor\alpha n\rfloor$.
    Let $p/q\in\QQ_{>0}$ where $p,q\in\ZZ_{>0}$ are coprime integers.
    We compute the values of $\lim_{\alpha\to\frac{p}{p+q}^+} c_{\alpha}(n)$ 
    at $n=-1$, $n=0$ and $n=p+q-1$:
    \begin{align*}
         \lim_{\alpha\to\frac{p}{p+q}^+} c_{\alpha}(-1) 
        &= \lim_{\alpha\to\frac{p}{p+q}^+} 
           \lfloor\alpha(-1+1)\rfloor-\lfloor\alpha (-1)\rfloor
        = 0 - (-1)
        = 1,\\
         \lim_{\alpha\to\frac{p}{p+q}^+} c_{\alpha}(0) 
        &= \lim_{\alpha\to\frac{p}{p+q}^+} 
           \lfloor\alpha(0+1)\rfloor-\lfloor\alpha (0)\rfloor
        = 0 - 0 = 0,\\
         \lim_{\alpha\to\frac{p}{p+q}^+} c_{\alpha}(p+q-1) 
        &= \lim_{\alpha\to\frac{p}{p+q}^+} 
           \lfloor\alpha(p+q-1+1)\rfloor-\lfloor\alpha (p+q-1)\rfloor
        = p - (p-1) = 1.
    \end{align*}
    For $0\leq n\leq p+q-1$,
    we have
$\lim_{\alpha\to\frac{p}{p+q}^+}c_{\alpha}(1)
    \cdots
    c_{\alpha}(p+q-1)=0m1\in\{0,1\}^*$ is the lower
    Christoffel word of slope $p/q$.  
    We now prove that $p+q$ is a period of 
    $n\mapsto\lim_{\alpha\to\frac{p}{p+q}^+} c_{\alpha}(n)$ 
    on the domain $\ZZ_{\geq0}$.
    Let $n\geq0$, we have
    \begin{align*}
        \lim_{\alpha\to\frac{p}{p+q}^+} c_{\alpha}(n+p+q) 
        &= \lim_{\alpha\to\frac{p}{p+q}^+} 
           \lfloor\alpha(n+p+q+1)\rfloor-\lfloor\alpha (n+p+q)\rfloor\\
        &= \lim_{\alpha\to\frac{p}{p+q}^+} 
           \lfloor\alpha(n+1)\rfloor+p-\lfloor\alpha n\rfloor-p\\
        &= \lim_{\alpha\to\frac{p}{p+q}^+} \lfloor\alpha(n+1)\rfloor-\lfloor\alpha n\rfloor
        = \lim_{\alpha\to\frac{p}{p+q}^+} c_{\alpha}(n).
    \end{align*}
    From \Cref{lem:symmetry-lower-to-upper}, we have
    $c_{\alpha}(n)=c_{\alpha}(-n-1)$ for every $n\in\ZZ\setminus\{-1,0\}$ and this shows the
    first equality since $m$ is a palindrome:
    \[
        \lim_{\alpha\to\frac{p}{p+q}^+} c_\alpha 
        = \leftidx{^\infty}(1m0)(1m1)\zeropoint(0m1)(0m1)^\infty.
    \]
    The other equalities are proved similarly.
\end{proof}

\begin{remark}
In general, the pair
    $\left(\leftidx{^\infty}(1m0) (1m1)\zeropoint(0m1) (0m1)^\infty,
    \leftidx{^\infty}(1m0) (1m0)\zeropoint(1m1) (0m1)^\infty\right)$
    is not indistinguishable.
For instance, when $m=0011$ the asymptotic pair
\begin{align*}
    x&= \leftidx{^\infty}(100110) (100111)\zeropoint(000111) (000111)^\infty,\\
    y&= \leftidx{^\infty}(100110) (100110)\zeropoint(100111) (000111)^\infty
\end{align*}
is not indistinguishable because the pattern $00111$ appears in $x$ intersecting the
difference set $F=\{-1,0\}$, but it does not appear in $y$ intersecting the
difference set $F$.
\end{remark}

    \subsection{Non-recurrent indistinguishable asymptotic pairs}

    We first prove that limits of pairs consisting of lower and upper characteristic Sturmian sequences whose slope tends towards a rational number are indistinguishable asymptotic pairs.
    Then, we show 
    in~\Cref{thm:sturmian_characterization_non_recurrent}
    that non-recurrent indistinguishable asymptotic pairs 
    whose difference set is of size 2
    are limits of Sturmian sequences.
    This result together with~\Cref{thm:sturmian_characterization}
    implies~\Cref{thm:sturmian_characterization_Z1}.

\begin{proposition}\label{prop:limit-of-sturmian-are-indistinguishable}
    Let $p,q\in\ZZ_{\geq0}$ be coprime integers.
    The limits of pairs consisting of lower and upper characteristic Sturmian sequences whose slope tends towards
    $p/(p+q)$ from above or from below
\[
    \left(\lim_{\alpha\to\frac{p}{p+q}^+} c_\alpha, 
          \lim_{\alpha\to\frac{p}{p+q}^+} c'_\alpha\right)
          \qquad
          \text{ and }
          \qquad
    \left(\lim_{\alpha\to\frac{p}{p+q}^-} c_\alpha, 
          \lim_{\alpha\to\frac{p}{p+q}^-} c'_\alpha\right).
\]
form two indistinguishable asymptotic pairs in $\ZZ$.
\end{proposition}

\begin{proof}
    From \Cref{lem:limit-of-sturmian-as-christoffel}, we observe that
    $\lim_{\alpha\to\frac{p}{p+q}^+} c_\alpha$
    and
    $\lim_{\alpha\to\frac{p}{p+q}^+} c'_\alpha$
    form an asymptotic pair whose difference set is $\{-1,0\}$.
    From \Cref{prop:limit-of-indist-is-indist},
    the property of being an indistinguishable pair is preserved by the limit.
    Therefore, it is an indistinguishable asymptotic pair.
    The same holds for the second pair.
\end{proof}

    \begin{theorem}\label{thm:sturmian_characterization_non_recurrent}
		Let $x,y\in\{\symb{0},\symb{1}\}^\ZZ$ and assume that $x$ is not recurrent.
		The pair $(x,y)$ is an indistinguishable asymptotic pair
		with difference set $F=\{-1,0\}$ such that $x_{-1}x_{0}=\symb{10}$
		and $y_{-1}y_{0}=\symb{01}$ if and only if
			there exist coprime nonnegative integers
			$p,q$ such that
			$(x,y)$ is the limit of pairs consisting of lower and upper characteristic Sturmian
			sequences whose slope tends to $p/(p+q)\in[0,1]\cap\QQ$ either
			\begin{itemize}
				\item from above, that is,
				$(x,y)=\lim_{\alpha\to\frac{p}{p+q}^+}(c_{\alpha},{c}'_{\alpha})$
				and $x=\sigma^{p+q}(y)$ is a shift of $y$, or,
				\item from below, that is,
				$(x,y)=\lim_{\alpha\to\frac{p}{p+q}^-}(c_{\alpha},{c}'_{\alpha})$
				and $y=\sigma^{p+q}(x)$ is a shift of $x$.
			\end{itemize}
	\end{theorem}

	\begin{proof}
        Let $p,q\in\ZZ_{\geq0}$ be coprime integers.
        From \Cref{prop:limit-of-sturmian-are-indistinguishable},
        the limits of pairs consisting of lower and upper characteristic Sturmian
        sequences whose slope tends to
        $p/(p+q)$ from above or from below form an indistinguishable asymptotic
        pair $(x,y)$ with $F=\{-1,0\}$ as a difference set where
        $x_{-1}x_{0}=\symb{10}$ and $y_{-1}y_{0}=\symb{01}$.

        Since $x$ is not recurrent,
        from \Cref{prop:recurrence}, we have that $x$ is a shift of $y$, i.e.
        $x=\sigma^k(y)$ for some $k\in\ZZ$. 
        We know that $k\neq0$ since $x\neq y$.
        If $k=1$, then
            $x = \leftidx{^\infty}01.00^\infty$
            and
            $y = \leftidx{^\infty}00.10^\infty$.
        From \Cref{lem:limit-of-sturmian-as-christoffel}, we conclude that
        \[
            x=\lim_{\alpha\to 0^+} c_\alpha 
            \qquad
            \text{ and }
            \qquad
            y=\lim_{\alpha\to 0^+}c'_\alpha
        \]
        are limits of pairs consisting of lower and upper characteristic Sturmian
        sequences whose slope tends to
        $0$ from above.
        Similarly, if $k=-1$, then
            $x = \leftidx{^\infty}11.01^\infty$ and
            $y = \leftidx{^\infty}10.11^\infty$.
        From \Cref{lem:limit-of-sturmian-as-christoffel}, we conclude that
            $x=\lim_{\alpha\to 1^-} c_\alpha $
            and
            $y=\lim_{\alpha\to 1^-}c'_\alpha$
        are limits of pairs consisting of lower and upper characteristic Sturmian
        sequences whose slope tends to
        $1$ from below.

        Assume now that $k\geq2$.
        Thus $y_{k-1}y_k=\symb{10}$.
        But $y_{k-1}y_k=x_{k-1}x_k$ so that $x_{k-1}x_k=\symb{10}$.
        Thus $y_{nk-1}y_{nk}=x_{nk-1}x_{nk}=\symb{10}$ for all $n>0$.
        Moreover $x_{-k-1}x_{-k}=y_{-1}y_{0}=\symb{01}$
        and
        $x_{nk-1}x_{nk}=y_{nk-1}y_{nk}=\symb{01}$ for all $n<0$.
        Let $m=x_1\cdots x_{k-1}$.
        We have $m=x_{nk+1}\cdots x_{nk-2}=y_{nk+1}\cdots y_{nk-2}$ for every $n\in\ZZ$.
        Thus we have
        \begin{align*}
            x &= \leftidx{^\infty}(1m0)(1m1)\zeropoint(0m1)(0m1)^\infty,\\
            y &= \leftidx{^\infty}(1m0)(1m0)\zeropoint(1m1)(0m1)^\infty.
        \end{align*}
        We observe that the factor $1m0$ appears in $y$ intersecting the difference set $F$.
        By the hypothesis, it must appear in $x$ intersecting the difference set $F$.
        Thus $1m0$ is a factor of $1m1.0m1$, but certainly not as a prefix.
        Therefore $1m0$ is a factor of $m1.0m1$.
        We conclude that $1m0$ is a factor of $0m1.0m1=(0m1)^2$.
        This implies that $1m0$ and $0m1$ are conjugate.
        From Pirillo's Theorem, we conclude that $0m1$ is a lower Christoffel word
        of slope $p/q$ for some coprime integers $p,q\in\ZZ_{\geq0}$ satisfying $p+q=k$.
        From \Cref{lem:limit-of-sturmian-as-christoffel}, we conclude that
        \[
            x=\lim_{\alpha\to\frac{p}{p+q}^+} c_\alpha 
            \qquad
            \text{ and }
            \qquad
            y=\lim_{\alpha\to\frac{p}{p+q}^+}c'_\alpha
        \]
        are limits of pairs consisting of lower and upper characteristic Sturmian
        sequences whose slope tends to
        $p/(p+q)$ from above.

        The proof for $k\leq-2$ follows the same line as when $k\geq2$ or can even
        be deduced from it by considering the reversal of $x$ and $y$.
        We obtain that
            $x=\lim_{\alpha\to\frac{p}{p+q}^-} c_\alpha $
            and
            $y=\lim_{\alpha\to\frac{p}{p+q}^-}c'_\alpha$
        are the limits of a sequence of lower and upper characteristic Sturmian sequences respectively whose slope converges to 
        $p/(p+q)$ from below.
	\end{proof}

    We may now deduce \Cref{thm:sturmian_characterization_Z1}.

    \begin{proof}[Proof of \Cref{thm:sturmian_characterization_Z1}]
        We have two cases to consider depending on whether $x$ is recurrent or not.
        If $x$ is recurrent, then from \Cref{thm:sturmian_characterization},
        we have that 
        the pair $(x,y)$ is an indistinguishable asymptotic pair with
        difference set $\{-1,0\}$ 
such that $x_{-1}x_{0}=\symb{10}$
		and $y_{-1}y_{0}=\symb{01}$ if and only if
            there exists $\alpha\in[0,1]\setminus\QQ$ such that
			$x={c}_{\alpha}$ and $y={c}'_{\alpha}$ are the lower and upper characteristic
			Sturmian sequences respectively.
        In this case, we consider the constant sequence $(\alpha_n)_{n\in\NN}$
        where $\alpha_n=\alpha$ for every $n\in\NN$.
        We have
            $x = c_\alpha=\lim_{n\to\infty} c_{\alpha_n}$
            and
            $y = c'_\alpha=\lim_{n\to\infty} c'_{\alpha_n}$.

        If $x$ is not recurrent, 
        then from \Cref{thm:sturmian_characterization_non_recurrent}
        the pair $(x,y)$ is an indistinguishable asymptotic pair with
        difference set $\{-1,0\}$ such that $x_{-1}x_{0}=\symb{10}$
		and $y_{-1}y_{0}=\symb{01}$ if and only if
        there exist coprime nonnegative integers $p$ and $q$ such that
		$(x,y)$ is the limit of a sequence of pairs of lower and upper characteristic Sturmian
        sequences whose slope tends toward the rational slope $p/(p+q)\in[0,1]\cap\QQ$ from above or
        from below.
        Let $(\alpha_n)_{n\in\NN}$ be the sequence defined as
        $\alpha_n=\frac{p}{p+q}+\frac{1}{\sqrt{2}n}$
        if the limit is from above or
        $\alpha_n=\frac{p}{p+q}-\frac{1}{\sqrt{2}n}$
        if the limit is from below.
        Then $\alpha_n \in [0,1)\setminus \QQ$ for all $n\in\NN$ and
        $(x,y)=\lim_{n\to\infty}(c_{\alpha_n},{c}'_{\alpha_n})$.
	\end{proof}

	\section{Indistinguishable asymptotic pairs on an arbitrary alphabet}\label{sec:resultsZ}
	
	The purpose of this section is to prove~\Cref{thm:characterization_Z} and hence provide a full characterization of indistinguishable asymptotic pairs in the case where the alphabet and difference set are arbitrary.~\Cref{thm:characterization_Z} will follow from~\Cref{prop:if_main_thm,prop:only_if_main_thm,prop:only_if_main_thm_non_recurrent}.

	\subsection{Substitutions preserve indistinguishability}
	
	We shall now show that indistinguishable asymptotic pairs are preserved under substitutions.

	\begin{definition}
		Let $\Sigma,\Gamma$ be alphabets. A \define{substitution} is a map $\varphi\colon \Sigma \to \Gamma^+$.
		
		The extension of $\varphi$ to a morphism from $\Sigma^+ \to \Gamma^+$ by concatenation is denoted (by abuse of notation) again $\varphi$. Moreover, every substitution induces a continuous map denoted (again, by abuse of notation) $\varphi \colon \Sigma^\ZZ \to \Gamma^\ZZ$ in the following way:
		\[ \varphi(x) \isdef \dots \varphi(x_{-5})\varphi(x_{-4})\varphi(x_{-3})\varphi(x_{-2})\varphi(x_{-1}) \zeropoint \varphi(x_{0})\varphi(x_{1})\varphi(x_{2})\varphi(x_{3})\varphi(x_{4})\dots  \]
	\end{definition}
	
	\begin{lemma}\label{lem:substitution_preserves_SI}
Let $\varphi\colon \Sigma \to \Gamma^+$ be a substitution and $x,y \in \Sigma^{\ZZ}$. 
        If $(x,y)$ is an indistinguishable asymptotic pair such that its difference set $F$ is contained in $\llbracket 0, k-1\rrbracket$, then
        $(\varphi(x),\varphi(y))$ is an indistinguishable asymptotic pair.
	\end{lemma}
	
	\begin{proof} From $F \subseteq \llbracket 0, k-1\rrbracket$ it follows immediately that for all $m <0$, $\varphi(x)_m = \varphi(y)_m$. Let $a \in \Sigma$. As $\Delta_a(x,y) = 0$, we deduce that $a$ appears the same number of times $N_a$ in both $x$ and $y$ in $\llbracket 0, k-1\rrbracket$.
	We deduce that  \[K \isdef |\varphi(x_{0})\dots \varphi(x_{k-1})| = \sum_{i = 0}^{k-1}|\varphi(x_i)| = \sum_{a \in \Sigma}N_a|\varphi(a)| =\sum_{i = 0}^{k-1}|\varphi(y_{i})| = |\varphi(y_{0})\dots \varphi(y_{k-1})|,\] and thus $\varphi(x)_m = \varphi(y)_m$ for every $m \geq K$.
	This shows that $\varphi(x)$ and $\varphi(y)$ are asymptotic and their difference set is contained in $D = \llbracket 0, K-1\rrbracket$.
	
	As $\Lcal(x) = \Lcal(y)$, we conclude that $\Lcal(\varphi(x)) = \Lcal(\varphi(y))$.	
		Fix $w \in \Lcal_n(\varphi(x))$.
		It suffices to show that 
        \begin{equation}\label{eq:occurrence-permutation}
            \#(\occ_{w}(\varphi(x)) \cap (D - \llbracket 0, n-1\rrbracket)) = \#(\occ_{w}(\varphi(y)) \cap (D - \llbracket 0, n-1\rrbracket)).
        \end{equation}

Every $i \in \occ_{w}(\varphi(x)) \cap (D- \llbracket 0,n-1\rrbracket)$ can be uniquely associated to a word $u_i \in \Sigma^+$ and a nonnegative integer $k_i$ such that $\sigma^{i-k_i}(\varphi(x)) \in [\varphi(u_i)]$, $\varphi(u_i)_{k_i}\dots \varphi(u_i)_{k_i+n-1}=w$ and so that $u_i$ is the shortest such word.
        As $\Delta_{u_i}(x,y)=0$, we have that $u_i$ occurs the same number of times in $x$ and $y$ in the support $F-\llbracket 0,|u_i|-1\rrbracket$. Therefore there is a bijection between $\occ_{u_i}(x) \cap (F-\llbracket 0,|u_i|-1\rrbracket)$ and $\occ_{u_i}(y) \cap (F-\llbracket 0,|u_i|-1\rrbracket)$ which induces a bijection between the set \[ A_{w,u,k}(x) =\{ i \in (\occ_{w}(\varphi(x)) \cap (D- \llbracket 0,n-1\rrbracket)  ) : i \mbox{ is associated to the pair } (u,k)\in \Sigma^+ \times \NN \},    \]
and the set
 \[ A_{w,u,k}(y) = \{i \in (\occ_{w}(\varphi(y)) \cap (D- \llbracket 0,n-1\rrbracket)  ) : i \mbox{ is associated to the pair } (u,k)\in \Sigma^+ \times \NN \}.    \]

As $\occ_{w}(\varphi(x)) \cap (D- \llbracket 0,n-1\rrbracket)$ can be written as the union of the $A_{w,u,k}(x)$ over all pairs $(u,k)$, and the same holds exchanging $x$ by $y$, we obtain that
\Cref{eq:occurrence-permutation} holds.
Thus $\Delta_w(\varphi(x),\varphi(y)) =0$. By~\Cref{prop:trivialite}, this implies that $\varphi(x),\varphi(y)$ form an indistinguishable asymptotic pair.
		\end{proof}
		
    We may now prove part of~\Cref{thm:characterization_Z} 
    based on~\Cref{prop:Sturmian_are_SI} and~\Cref{lem:substitution_preserves_SI}.
	
	\begin{proposition}\label{prop:if_main_thm}
		Let $\alpha$ be irrational and ${c}_{\alpha},{c}'_{\alpha}$ be the lower and upper characteristic words of slope $\alpha$ respectively. For any substitution $\varphi\colon \{\symb{0},\symb{1}\} \to \Sigma^+$ the sequences
		$\varphi(\sigma^{1}({c}_{\alpha}))$ and $\varphi(\sigma^{1}({c}'_{\alpha}))$ form an indistinguishable asymptotic pair.
	\end{proposition}
	
	\begin{proof} By~\Cref{prop:Sturmian_are_SI}, we have that ${c}_{\alpha},{c}'_{\alpha}$ form a non-trivial indistinguishable asymptotic pair. By~\Cref{prop:shifted_SI-for-Z}, $\sigma^{-1}({c}_{\alpha})$ and $\sigma^{-1}({c'}_{\alpha})$ are also a non-trivial indistinguishable asymptotic pair with difference set $F = \{0,1\}$. By~\Cref{lem:substitution_preserves_SI}, we have that $\varphi(\sigma^{-1}({c}_{\alpha})),\varphi(\sigma^{-1}({c}'_{\alpha}))$ is an indistinguishable asymptotic pair.
		
	Note that if we let $m = |\varphi(\symb{0})|+|\varphi(\symb{1})|$ then $\sigma^m(\varphi(\sigma^{-1}({c}_{\alpha}))) = \varphi(\sigma^{1}({c}_{\alpha}))$ and $\sigma^m(\varphi(\sigma^{-1}({c}'_{\alpha}))) = \varphi(\sigma^{1}({c}'_{\alpha}))$. Then again by~\Cref{prop:shifted_SI-for-Z}, we obtain that $\varphi(\sigma^{1}({c}_{\alpha}))$ and $\varphi(\sigma^{1}({c}'_{\alpha}))$ form an indistinguishable asymptotic pair.\end{proof}
	
	Note that in~\Cref{prop:if_main_thm} we do not ensure that after applying the substitution the words remain non-trivial.
	For instance, we may consider a substitution sending all symbols to a fixed symbol and thus trivialize the pair.

    \subsection{Derived sequences preserve indistinguishability}\label{sec:derived-sequences}

     We shall find a sequence of inverse substitutions which will allow us to ``desubstitute'' the asymptotic pair until we arrive to a Sturmian sequence. The main tool is the notion of derived sequence introduced by Durand~\cite{MR1489074}.
	
	\begin{definition}\label{def:derived_seq}
		Let $x \in \Sigma^{\ZZ}$ and $w \in \Lcal(x)$ which appears with bounded gaps in $x$. Let $\left \{i_{k} \right\}_{k \in \ZZ}$ be the enumeration of $\occ_w(x)$ which is strictly increasing and such that $i_0$ is the smallest value of $\occ_w(x)$ such that $i_0 > -|w|$
		
		The \define{derived sequence} $D_{w}(x) \in (\texttt{RW}_w(x))^{\ZZ}$ is the sequence given by \[ (D_w(x))_k = x_{i_{k}}\dots x_{i_{k+1}-1}. \]
	\end{definition}
	
	The derived sequence of a uniformly recurrent sequence is also uniformly recurrent. Note that the alphabet of the derived sequence consists of symbols in $\texttt{RW}_w(x)$ which formally are words in $\Sigma^*$. It is possible to recover the original sequence $x$ (up to a shift) by applying the morphism $\varphi : \texttt{RW}_w(x) \to \Sigma^*$ such that $\varphi(u) = u_0\dots u_{|u|-1}$.
	
	\begin{lemma}\label{lem:derived_SI}
		Let $x,y \in \Sigma^*$ 
        and assume that $a\in\Sigma$ appears with bounded gaps in $x$.
        If $(x,y)$ is an indistinguishable asymptotic pair 
        whose difference set $F$ is contained in $\llbracket 0, k-1\rrbracket$,
        then $(D_a(x)$,$D_a(y))$ is an indistinguishable asymptotic pair.
		
        Moreover, the difference set of $(D_a(x)$,$D_a(y))$ is contained in $\llbracket 0,N_a\rrbracket$ where \[N_a = \#(\{i \in \llbracket 0, k-1\rrbracket : x_i = a \}).\]
	\end{lemma}
	
	\begin{proof}
		Rewrite the sets $\occ_a(x)$ and $\occ_a(y)$ in increasing order as sequences $\left \{i_{t} \right \}_{t \in \ZZ}$ and $\left \{j_{t} \right \}_{t \in \ZZ}$ as in~\Cref{def:derived_seq} respectively.  As $F \subseteq \llbracket 0,k-1\rrbracket$, and $x,y$ are asymptotic, we have that $i_{t} = j_{t}$ for every $t < 0$ and so $(D_a(x))_t = (D_a(y))_t$ for every $t < 0$. 
		
		As $\Delta_a(x,y)=0$, $a$ occurs the same number of times $N_{a}$ in the interval $\llbracket 0,k-1\rrbracket$ in $x$ and $y$, Therefore using again that $x,y$ are asymptotic, we have that $i_{t} = j_{t}$ for every $t \geq N_{a}+1$ and thus $(D_a(x))_t = (D_a(y))_t$ for every $t \geq N_{a}+1$. This shows that $D_a(x)$ and $D_a(y)$ are asymptotic and that their difference set is contained in $\llbracket0,N_a\rrbracket$.
		
		Let $\varphi : \texttt{RW}_a(x) \to \Sigma^*$ be the morphism such that $\varphi(u) = u_0\dots u_{|u|-1}$.
		Given a word $w = w_1\dots w_m \in (\texttt{RW}_a(x))^*$, let $|\varphi(w)| = \sum_{i=1}^m |\varphi(w_i)|$. It follows that
		
		\begin{align*}
		\Delta_{w}(D_a(x),D_a(y)) & = \sum_{\ell=-(|w|-1)}^{N_a}\indicator{[w]}(\sigma^{\ell}(D_a(y)))-\indicator{[w]}(\sigma^{\ell}(D_a(x)))\\
		&  = \sum_{ \ell = -(|\varphi(w)|-1)}^{k-1} \indicator{[\varphi(w)]}(\sigma^{\ell}(y)) - \indicator{[\varphi(w)]}(\sigma^{\ell}(x))  \\
		&  = \Delta_{\varphi(w)}(x,y) = 0.
		\end{align*}
		
		It follows that $D_a(x)$,$D_a(y)$ form an indistinguishable asymptotic pair.\end{proof}
	
	\begin{remark}
		An analogous statement holds if instead of considering $a \in \Sigma$ we take an arbitrary $w \in \Lcal(x)$ and we consider the pair $D_w(x),D_w(y)$. We shall not need this general statement.
	\end{remark}
	
	\subsection{Proof of~\Cref{thm:characterization_Z}}
	We will first show that, as long as the difference set of an indistinguishable asymptotic pair is contained in an interval of length at least $3$, we can use derived sequences to construct a new indistinguishable pair whose difference set is contained in a strictly smaller interval. This will later provide a way to reduce the general case to the case where the difference set is $\{-1,0\}$.
	
	\begin{lemma} \label{lem:iteration}
		Suppose $x,y \in \Sigma^{\ZZ}$ is a non-trivial indistinguishable asymptotic pair whose difference set is contained in an interval $F = \llbracket 0, k-1\rrbracket$. 
        If $x$ is recurrent, there is $a \in \Sigma$ such that $D_{a}(x)$ and
        $D_{a}(y)$ form a non-trivial indistinguishable asymptotic pair with a
        difference set contained in the interval $\llbracket 0, \lfloor\frac{k}{2}
        \rfloor\rrbracket$.
	\end{lemma}

	\begin{proof}
		As $x,y$ are non-trivial, we have that $\#(\Sigma)\geq 2$. Let $a \in \Sigma$ be the symbol such that $\occ_a(x) \cap\llbracket 0, k-1\rrbracket$ is the smallest. By the pigeonhole principle, $N_a \isdef \#(\occ_a(x) \cap\llbracket 0, k-1\rrbracket)\leq \lfloor\frac{k}{2} \rfloor$.
		
		By~\Cref{lem:rec_implies_urec}, both $x$ and $y$ are uniformly recurrent and so the sequences $D_{a}(x)$ and $D_{a}(y)$ are well defined. By~\Cref{lem:derived_SI} $D_{a}(x)$ and $D_{a}(y)$ form an indistinguishable asymptotic pair with difference set contained in $\llbracket 0,N_a\rrbracket \subseteq \llbracket 0, \lfloor\frac{k}{2}
		\rfloor\rrbracket$, which is clearly non-trivial as $x,y$ are non-trivial.
	\end{proof}

    	\begin{proposition}\label{prop:only_if_main_thm}
        Let $x,y \in \Sigma^{\ZZ}$ be a non-trivial
        indistinguishable asymptotic pair. If $x$ is recurrent, then there exists 
        $\alpha \in [0,1) \setminus \QQ$, 
        a substitution $\varphi\colon \{\symb{0},\symb{1}\} \to \Sigma^+$ 
        and an integer $m \in \ZZ$ such that 
        \[ 
        \{x,y\} = \{\sigma^m(\varphi(\sigma(c_{\alpha}))),\sigma^m(\varphi(\sigma(c'_{\alpha})))  \} 
        \]
        where $c_{\alpha}$ and $c'_{\alpha}$ are the lower and upper
        characteristic Sturmian sequences of slope $\alpha$ respectively.
	\end{proposition}		
	
	\begin{proof}
        Let $\llbracket\ell,\ell+k-1\rrbracket$ be the smallest interval
        containing the difference set $F$ of $(x,y)$.
        Since $(x,y)$ is non-trivial and indistinguishable, we have that $k \geq 2$. We shall proceed by induction on
        $k$.
		If $k = 2$, by~\Cref{prop:complexity-bound} we have that the alphabet has size at most $\# \Lcal_1(x) \leq 1+k-1 = 2$. Therefore, up to a relabeling of the alphabet and a shift, we have an asymptotic pair of sequences which satisfies the assumptions of~\Cref{thm:sturmian_characterization} and therefore
        $\{\sigma^{\ell+1}(x),\sigma^{\ell+1}(y)\} = \{\varphi(c_{\alpha}),\varphi(c'_{\alpha})\}$ 
        are the lower and upper characteristic Sturmian sequences
        for some $\alpha \in [0,1) \setminus \QQ$
        up to some function $\varphi \colon \{\symb{0},\symb{1}\} \to \Sigma$.
        We conclude that $\{x,y\} = \{\sigma^{-\ell-2}(\varphi(\sigma(c_{\alpha}))),\sigma^{-\ell-2}(\varphi(\sigma(c'_{\alpha})))\}$.
		
		Now suppose $k \geq 3$ and the result holds for all $2 \leq j < k$. By~\Cref{prop:shifted_SI-for-Z}, $x^* = \sigma^{\ell}(x)$ and $y^* = \sigma^{\ell}(y)$ are an indistinguishable asymptotic pair whose difference set is contained in $\llbracket 0,k-1\rrbracket$. By~\Cref{lem:iteration} there is $a \in \Sigma$ such that $x' \isdef D_{a}(x^*)$ and $y' \isdef D_{a}(y^*)$ are a non-trivial indistinguishable pair in $(\texttt{RW}_a(x))^{\ZZ}$ and their difference set is contained in the interval $\llbracket 0,\lfloor\frac{k}{2}\rfloor \rrbracket$. As $k \geq 3$, we have that $\lfloor\frac{k}{2}\rfloor < k-1$ and thus the result holds for $x',y'$. It follows that there is a substitution $\varphi' \colon \{\symb{0},\symb{1}\} \to (\texttt{RW}_a(x) )^+$ and $m' \in \ZZ$ so that 
		\[ \{ x',y'\} = \{\sigma^{m'}(\varphi'(\sigma(c_{\alpha}))), \sigma^{m'}(\varphi'(\sigma(c'_{\alpha})))  \}.   \]
		
		Let $\phi \colon \texttt{RW}_a(x) \to \Sigma^+$ and $s \in \ZZ$ be respectively the substitution and integer such that $\sigma^{s}(\phi(x'))=x^*$ and $\sigma^{s}(\phi(y'))=y^*$. Let $\varphi \isdef \phi \circ \varphi'$.
		Note that the difference set of $x',y'$ is contained in $\llbracket 0,\lfloor\frac{k}{2}\rfloor \rrbracket$ and the difference set of $\sigma(\varphi'(c_{\alpha})),\sigma(\varphi'(c'_{\alpha}))$ is contained in $\llbracket-(|\varphi'(\symb{0})|+|\varphi'(\symb{1})|),-1 \rrbracket$. Let us first show that there is $N \in \NN$ so that $\phi(\sigma^{m'}(\varphi'(\sigma(c_{\alpha})))) = \sigma^{-N}(\varphi(\sigma(c_{\alpha})))$ and $\phi(\sigma^{m'}(\varphi'(\sigma(c'_{\alpha})))) = \sigma^{-N}(\varphi(\sigma(c'_{\alpha})))$.

		Let $K \in \NN$ be the smallest positive integer such that $\llbracket-K,-1\rrbracket$ contains the difference set of $\varphi'(\sigma(c_{\alpha})),\varphi'(\sigma(c'_{\alpha}))$. As the difference set of $x',y'$ is contained in $\llbracket0,\lfloor\frac{k}{2}\rfloor \rrbracket$, we obtain that $m' \leq -K+1$. Consider the words \begin{align*}
		w_1 & = \varphi'(\sigma(c_{\alpha}))_{m'}\dots\varphi'(\sigma(c_{\alpha}))_{-1}\\
		w_2 & =  \varphi'(\sigma(c'_{\alpha}))_{m'}\dots\varphi'(\sigma(c'_{\alpha}))_{-1}.
		\end{align*}
		By construction, every symbol in $\texttt{RW}_a(x)$ occurs the same number of times in $w_1$ and $w_2$. Letting $N = |\phi(w_1)| = |\phi(w_2)|$ we obtain that $\phi(\sigma^{m'}(\varphi'(\sigma(c_{\alpha})))) = \sigma^{-N}(\varphi(\sigma(c_{\alpha})))$ and $\phi(\sigma^{m'}(\varphi'(\sigma(c'_{\alpha})))) = \sigma^{-N}(\varphi(\sigma(c'_{\alpha})))$.
		
        Finally, we conclude that  \[\{x,y\} = \{ \sigma^{\ell+s}(\phi(x')), \sigma^{\ell+s}(\phi(y'))\} =  \{\sigma^{\ell+s-N}(\varphi(\sigma(c_{\alpha}))),\sigma^{\ell+s-N}(\varphi(\sigma(c'_{\alpha})))\}\]
        which is what we wanted to prove.
    \end{proof}

    We deal with the case when $x$ is non-recurrent in the following proposition.
    
    \begin{proposition} \label{prop:only_if_main_thm_non_recurrent}
    Let $x,y \in \Sigma^{\ZZ}$ be a non-trivial indistinguishable asymptotic pair.
	If $x$ is not recurrent, then there exists a substitution 
			$\varphi \colon \{\symb{0},\symb{1}\} \to \Sigma^+$  and an integer $m\in\ZZ$
			such that 
			\[\{x,y  \} = \{\sigma^m\varphi(\leftidx^{\infty}\symb{0}.\symb{1}\symb{0}^{\infty}), \sigma^m\varphi(\leftidx^{\infty}\symb{0}.\symb{0}\symb{1}\symb{0}^{\infty})  \}.\]
    \end{proposition}
    \begin{proof}
    By \Cref{prop:recurrence}, $x$ and $y$ lie on the same orbit, i.e., there exists $s\in \ZZ\setminus \{0\}$ with $x = \sigma^s(y)$. Possibly exchanging $x$ and $y$, we may assume $s > 0$.
    Let $r = \min\{i \in \ZZ : x_{i} \neq y_{i}\}$, then the difference set of $\sigma^{r}(x),\sigma^{r}(y)$ is contained in the interval $\llbracket 0,k-1\rrbracket$ for some $k > s$. Let us denote $x' = \sigma^{r}(x)$ and $y' = \sigma^{r}(y)$.
    
    A word $u$ that cannot be written as a repeated concatenation of another word is called \define{primitive}. Let $u$ be a primitive word such that $u^n= x'_{k}\dots x'_{k+s-1}$ for some positive integer $n$. Since $x'|_{\llbracket k,\infty \rrbracket} = y'|_{\llbracket k,\infty \rrbracket} = \sigma^s(x')|_{\llbracket k,\infty \rrbracket}$, we obtain that $x'$ and $y'$ are eventually periodic to the right, more precisely
    \[
    x'|_{\llbracket k-s,\infty \rrbracket} =  u^{\infty} \quad \mbox{ and }\quad y'|_{\llbracket k,\infty \rrbracket} =  u^{\infty}
    \]
    
    Similarly, let $w$ be a primitive word such that $w^{n'}= x'_{-s}\dots x'_{-1}$ for some positive integer $n'$. Since we have
    $ y'|_{ \llbracket -\infty, -1\rrbracket}= x'|_{ \llbracket -\infty, -1\rrbracket} = \sigma^{-s}(y')|_{ \llbracket -\infty, -1\rrbracket}$
    we obtain that $x'$ and $y'$ are eventually periodic to the left, more precisely
    \[
    x'|_{ \llbracket -\infty, -1\rrbracket} = \leftidx^{\infty}w \quad \mbox{ and } \quad  y'|_{ \llbracket -\infty, s-1\rrbracket} = \leftidx^{\infty}w
    \]
    
    As we took $k > s$, letting $v = x'_0 \dots x'_{k-s-1}$ we can write \[x' = \leftidx^{\infty}w\zeropoint vuu^{\infty}, \quad \mbox{ and } \quad y' = \leftidx^{\infty}w\zeropoint wvu^{\infty}.\] 
    We claim that $w$ and $u$ are conjugate. Indeed, as $x'$ and $y'$ are indistinguishable, we have that \[\#(\occ_w(x')\cap \llbracket-|w|+1,k-1\rrbracket) = \#(\occ_w(y')\cap \llbracket-|w|+1,k-1\rrbracket).   \]
    
    First, note that there cannot be any occurrence of $w$ in $y'$ on the indexes $\llbracket-|w|+1,-1\rrbracket$, otherwise we would have that $w$ occurs as a factor of $ww$ neither as a prefix nor suffix, which would contradict the primitivity of $w$. Also note that every occurrence of $w$ in $y'$ on an index $j \in \llbracket1,k-1\rrbracket$ can be mapped uniquely to an occurrence of $w$ in $x'$ on index $j-|w|$. 
    Therefore, the only remaining occurrence of $w$ in $y'$ at index $0$ must necessarily be mapped to an occurrence of $w$ in $x'$ as a factor of $uu$. Therefore $w$ is a factor of $uu$. 
 Since $u$ and $w$ have the same length, it
    implies that they are conjugate, that is, $u = pz$ and $w = zp$ for some words $z,p$. 
    Letting $t = vp$ we can rewrite $x'$ and $y'$ in the following way
     \[x' = \leftidx^{\infty}w\zeropoint tw^{\infty}, \quad \mbox{ and } \quad y' = \leftidx^{\infty}w\zeropoint wtw^{\infty}.\]
    Setting $\varphi: \symb{0} \mapsto w, \symb{1} \mapsto t$ we get $x' =\varphi(\leftidx^{\infty}\symb{0}.\symb{1}\symb{0}^{\infty})$ and  $y' = \varphi(\leftidx^{\infty}\symb{0}.\symb{0}\symb{1}\symb{0}^{\infty})$.
    Letting $m = -r$ we obtain \[ x = \sigma^m(\varphi(\leftidx^{\infty}\symb{0}.\symb{1}\symb{0}^{\infty})) \quad \mbox{ and } \quad y = \sigma^m(\leftidx^{\infty}\symb{0}.\symb{0}\symb{1}\symb{0}^{\infty}).  \]
    Which is what we wanted, modulo exchanging $x$ and $y$.
    \end{proof}

\begin{proof}[Proof of \Cref{thm:characterization_Z}]
    Let $\varphi\colon \{\symb{0},\symb{1}\}\to \Sigma^+$ be a substitution. From \Cref{prop:if_main_thm}, if $\alpha$ is irrational, then
    the sequences $\varphi(\sigma({c}_{\alpha}))$ and 
    $\varphi(\sigma({c}'_{\alpha}))$
    form an indistinguishable asymptotic pair and thus by ~\Cref{prop:shifted_SI-for-Z},
    the asymptotic pair
    $\sigma^m\varphi(\sigma({c}_{\alpha})), 
    \sigma^m\varphi(\sigma({c}'_{\alpha}))$
    is indistinguishable. Similarly, as $\leftidx^{\infty}\symb{0}.\symb{1}\symb{0}^{\infty},\leftidx^{\infty}\symb{0}.\symb{0}\symb{1}\symb{0}^{\infty}$ is indistinguishable, we have that for every integer $m\in\ZZ$, the pair
    $\sigma^m\varphi(\leftidx^{\infty}\symb{0}.\symb{1}\symb{0}^{\infty}),
    \sigma^m\varphi(\leftidx^{\infty}\symb{0}.\symb{0}\symb{1}\symb{0}^{\infty})$
    is indistinguishable.

    Conversely, if $x$ is recurrent, the result is proved in 
    \Cref{prop:only_if_main_thm}.
    If $x$ is non-recurrent, the result is proved in
    \Cref{prop:only_if_main_thm_non_recurrent}.
\end{proof}

\Addresses

\bibliographystyle{abbrv}
\bibliography{biblio}

\end{document}